\documentclass[12pt]{amsart}
\usepackage{amssymb}
\usepackage{amsmath}
\usepackage{amscd}
\usepackage[all]{xypic}

\allowdisplaybreaks[1]

\newtheorem{thm}{Theorem}[section]
\newtheorem{lem}[thm]{Lemma}
\newtheorem{cor}[thm]{Corollary}
\newtheorem{prop}[thm]{Proposition}

\newtheorem{tadothm}{Theorem}[section]
\newtheorem{tadotadothm}{Theorem}[section]
\newtheorem{tanakathm}{Theorem}[section]

\theoremstyle{definition}
\newtheorem{rem}[thm]{Remark}
\newtheorem{defn}[thm]{Definition}

\theoremstyle{remark}

\numberwithin{equation}{section}

\def\R{{\mathbb R}}
\def\Z{{\mathbb Z}}
\def\C{{\mathbb C}}

\def\Ind{\text{\rm Ind}}
\def\Coind{\text{\rm Coind}}
\def\Hom{\text{\rm Hom}}
\def\id{\text{\rm id}}

\begin{document}
\title[Harmonic Volumes of Hyperelliptic Curves]
{The Harmonic Volumes of Hyperelliptic Curves}

\author[Yuuki Tadokoro]{Yuuki Tadokoro}
\address{Department of Mathematical Sciences,
University of Tokyo, 3-8-1 Komaba, Meguro, Tokyo 153-8914, Japan}
\email{tado\char`\@ms.u-tokyo.ac.jp}


\begin{abstract}
We determine the harmonic volumes for all the hyperelliptic curves.
This gives a geometric interpretation
of a theorem established by A. Tanaka \cite{T}.
\end{abstract}

\maketitle
\section{Introduction}
Let $X$ be a compact Riemann surface of genus $g\ge 3$.
A harmonic volume $I$ of $X$ was introduced by 
B. Harris \cite{H-1}, using Chen's iterated integrals \cite{C}.
The aim of this paper is to determine the harmonic volumes of all the 
hyperelliptic curves, which are $2$-fold branched coverings of $\C P^1$.
As was already pointed out by Harris, some important algebraic cycles 
in the Jacobian variety $J(X)$ are related to $2I$, 
which vanishes for all the hyperelliptic curves.
The harmonic volumes of hyperelliptic curves, however,
 have been still unknown.
First of all, we give the statement of the main theorem of this paper.
 We denote by $H$ the first integral homology group of $X$.
Harris defined the harmonic volume $I$
 as a homomorphism $(H^{\otimes 3})^\prime\to \R/{\Z}$.
Here $(H^{\otimes 3})^\prime$ is a certain subgroup of $H^{\otimes 3}$.
See Section $2$ for the definition of $(H^{\otimes 3})^\prime$.
We denote by $C$ a hyperelliptic curve.
\begin{tadothm}
For any hyperelliptic curve $C$, let $\{x_i, y_i\}_{i=1,2, \ldots, g}$
 be a symplectic basis of $H=H_1(C; \Z)$ in Figure \ref{figure},
where $\iota$ is the hyperelliptic involution.
We denote by $z_i$ either $x_i$ or $y_i$.
Then,
$$
\mbox{\phantom{gi}}I(z_i\otimes z_j \otimes z_k) =
\left.
  \begin{array}{cc}
  0 & \mbox{for}\ i\neq j\neq k\neq i,
  \end{array}
\right.
$$
$$
I(x_i\otimes y_i\otimes z_k -x_{k+1}\otimes y_{k+1}\otimes z_k)=
\left\{
  \begin{array}{cc}
  {\displaystyle {1/{2}}} & \mbox{for}\ i<k, k=2, 3, \ldots, g-1\ \mbox{and}\ z_k=y_k,\\
  0 & \mbox{for}\ i\geq k+2, k=1, k=g\ \mbox{or}\ z_k=x_k.
  \end{array}
  \right.
$$
The elements $z_i\otimes z_j \otimes z_k$ and
 $x_i\otimes y_i\otimes z_k -x_{k+1}\otimes y_{k+1}\otimes z_k$
are the parts of a basis of $(H^{\otimes 3})^\prime$ whose harmonic
volumes depend on the complex structure of Riemann surfaces.
\begin{center}
\begin{figure}[htbp]
\unitlength 0.1in
\begin{picture}(56.00,16.00)(6.80,-22.00)
%
\special{pn 20}%
\special{pa 1800 600}%
\special{pa 5000 600}%
\special{fp}%
\special{pa 5000 2200}%
\special{pa 1800 2200}%
\special{fp}%
%
\special{pn 20}%
\special{ar 1800 1400 800 800  1.5707963 4.7123890}%
%
\special{pn 20}%
\special{ar 5000 1400 800 800  4.7123890 6.2831853}%
\special{ar 5000 1400 800 800  0.0000000 1.5707963}%
%
\special{pn 20}%
\special{ar 1960 1360 280 160  6.2831853 6.2831853}%
\special{ar 1960 1360 280 160  0.0000000 3.1415927}%
%
\special{pn 20}%
\special{ar 1960 1400 248 104  3.1415927 6.2831853}%
%
\special{pn 20}%
\special{ar 3400 1360 280 160  6.2831853 6.2831853}%
\special{ar 3400 1360 280 160  0.0000000 3.1415927}%
%
\special{pn 20}%
\special{ar 4840 1360 280 160  6.2831853 6.2831853}%
\special{ar 4840 1360 280 160  0.0000000 3.1415927}%
%
\special{pn 20}%
\special{ar 3400 1400 248 104  3.1415927 6.2831853}%
%
\special{pn 20}%
\special{ar 4840 1400 248 104  3.1415927 6.2831853}%
%
\special{pn 8}%
\special{sh 0.600}%
\special{ar 1040 1400 24 24  0.0000000 6.2831853}%
%
\special{pn 8}%
\special{sh 0.600}%
\special{ar 5760 1400 24 24  0.0000000 6.2831853}%
%
\special{pn 8}%
\special{sh 0.600}%
\special{ar 5144 1400 24 24  0.0000000 6.2831853}%
%
\special{pn 8}%
\special{sh 0.600}%
\special{ar 3704 1400 24 24  0.0000000 6.2831853}%
%
\special{pn 8}%
\special{sh 0.600}%
\special{ar 2264 1400 24 24  0.0000000 6.2831853}%
%
\special{pn 8}%
\special{sh 0.600}%
\special{ar 1656 1400 24 24  0.0000000 6.2831853}%
%
\special{pn 8}%
\special{sh 0.600}%
\special{ar 3096 1400 24 24  0.0000000 6.2831853}%
%
\special{pn 8}%
\special{sh 0.600}%
\special{ar 4536 1400 24 24  0.0000000 6.2831853}%
%
\special{pn 20}%
\special{sh 1}%
\special{ar 2680 1400 10 10 0  6.28318530717959E+0000}%
\special{sh 1}%
\special{ar 2600 1400 10 10 0  6.28318530717959E+0000}%
\special{sh 1}%
\special{ar 2760 1400 10 10 0  6.28318530717959E+0000}%
%
\special{pn 20}%
\special{sh 1}%
\special{ar 4120 1400 10 10 0  6.28318530717959E+0000}%
\special{sh 1}%
\special{ar 4040 1400 10 10 0  6.28318530717959E+0000}%
\special{sh 1}%
\special{ar 4200 1400 10 10 0  6.28318530717959E+0000}%
%
\special{pn 8}%
\special{pa 1000 1400}%
\special{pa 680 1400}%
\special{fp}%
\special{pa 5800 1400}%
\special{pa 6280 1400}%
\special{fp}%
%
\special{pn 8}%
\special{ar 6040 1400 160 480  3.7310362 6.2831853}%
\special{ar 6040 1400 160 480  0.0000000 2.3561945}%
%
\special{pn 8}%
\special{pa 5928 1720}%
\special{pa 5912 1680}%
\special{fp}%
\special{sh 1}%
\special{pa 5912 1680}%
\special{pa 5918 1749}%
\special{pa 5932 1730}%
\special{pa 5955 1734}%
\special{pa 5912 1680}%
\special{fp}%
\put(58.4000,-8.8000){\makebox(0,0)[lb]{$180^\circ$}}%
\put(60.4000,-20.4000){\makebox(0,0)[lb]{$\iota$}}%
\put(58.0000,-22.0000){\makebox(0,0)[lb]{$C$}}%
%
\special{pn 13}%
\special{ar 3400 1400 500 370  0.0000000 6.2831853}%
%
\special{pn 13}%
\special{ar 4850 1400 500 370  0.0000000 6.2831853}%
%
\special{pn 13}%
\special{ar 1960 1400 500 370  0.0000000 6.2831853}%
%
\special{pn 13}%
\special{ar 4850 1860 180 320  1.5707963 4.7123890}%
%
\special{pn 13}%
\special{ar 3400 1860 180 320  4.7123890 4.7603890}%
\special{ar 3400 1860 180 320  4.9043890 4.9523890}%
\special{ar 3400 1860 180 320  5.0963890 5.1443890}%
\special{ar 3400 1860 180 320  5.2883890 5.3363890}%
\special{ar 3400 1860 180 320  5.4803890 5.5283890}%
\special{ar 3400 1860 180 320  5.6723890 5.7203890}%
\special{ar 3400 1860 180 320  5.8643890 5.9123890}%
\special{ar 3400 1860 180 320  6.0563890 6.1043890}%
\special{ar 3400 1860 180 320  6.2483890 6.2963890}%
\special{ar 3400 1860 180 320  6.4403890 6.4883890}%
\special{ar 3400 1860 180 320  6.6323890 6.6803890}%
\special{ar 3400 1860 180 320  6.8243890 6.8723890}%
\special{ar 3400 1860 180 320  7.0163890 7.0643890}%
\special{ar 3400 1860 180 320  7.2083890 7.2563890}%
\special{ar 3400 1860 180 320  7.4003890 7.4483890}%
\special{ar 3400 1860 180 320  7.5923890 7.6403890}%
\special{ar 3400 1860 180 320  7.7843890 7.8323890}%
%
\special{pn 13}%
\special{ar 3400 1860 180 320  1.5707963 4.7123890}%
%
\special{pn 13}%
\special{ar 4850 1860 180 320  4.7123890 4.7603890}%
\special{ar 4850 1860 180 320  4.9043890 4.9523890}%
\special{ar 4850 1860 180 320  5.0963890 5.1443890}%
\special{ar 4850 1860 180 320  5.2883890 5.3363890}%
\special{ar 4850 1860 180 320  5.4803890 5.5283890}%
\special{ar 4850 1860 180 320  5.6723890 5.7203890}%
\special{ar 4850 1860 180 320  5.8643890 5.9123890}%
\special{ar 4850 1860 180 320  6.0563890 6.1043890}%
\special{ar 4850 1860 180 320  6.2483890 6.2963890}%
\special{ar 4850 1860 180 320  6.4403890 6.4883890}%
\special{ar 4850 1860 180 320  6.6323890 6.6803890}%
\special{ar 4850 1860 180 320  6.8243890 6.8723890}%
\special{ar 4850 1860 180 320  7.0163890 7.0643890}%
\special{ar 4850 1860 180 320  7.2083890 7.2563890}%
\special{ar 4850 1860 180 320  7.4003890 7.4483890}%
\special{ar 4850 1860 180 320  7.5923890 7.6403890}%
\special{ar 4850 1860 180 320  7.7843890 7.8323890}%
%
\special{pn 13}%
\special{ar 1950 1860 180 320  1.5707963 4.7123890}%
%
\special{pn 13}%
\special{ar 1950 1860 180 320  4.7123890 4.7603890}%
\special{ar 1950 1860 180 320  4.9043890 4.9523890}%
\special{ar 1950 1860 180 320  5.0963890 5.1443890}%
\special{ar 1950 1860 180 320  5.2883890 5.3363890}%
\special{ar 1950 1860 180 320  5.4803890 5.5283890}%
\special{ar 1950 1860 180 320  5.6723890 5.7203890}%
\special{ar 1950 1860 180 320  5.8643890 5.9123890}%
\special{ar 1950 1860 180 320  6.0563890 6.1043890}%
\special{ar 1950 1860 180 320  6.2483890 6.2963890}%
\special{ar 1950 1860 180 320  6.4403890 6.4883890}%
\special{ar 1950 1860 180 320  6.6323890 6.6803890}%
\special{ar 1950 1860 180 320  6.8243890 6.8723890}%
\special{ar 1950 1860 180 320  7.0163890 7.0643890}%
\special{ar 1950 1860 180 320  7.2083890 7.2563890}%
\special{ar 1950 1860 180 320  7.4003890 7.4483890}%
\special{ar 1950 1860 180 320  7.5923890 7.6403890}%
\special{ar 1950 1860 180 320  7.7843890 7.8323890}%
%
\special{pn 20}%
\special{pa 1770 1860}%
\special{pa 1770 1850}%
\special{fp}%
\special{sh 1}%
\special{pa 1770 1850}%
\special{pa 1750 1917}%
\special{pa 1770 1903}%
\special{pa 1790 1917}%
\special{pa 1770 1850}%
\special{fp}%
\special{pa 3220 1860}%
\special{pa 3220 1850}%
\special{fp}%
\special{sh 1}%
\special{pa 3220 1850}%
\special{pa 3200 1917}%
\special{pa 3220 1903}%
\special{pa 3240 1917}%
\special{pa 3220 1850}%
\special{fp}%
\special{pa 4670 1860}%
\special{pa 4670 1850}%
\special{fp}%
\special{sh 1}%
\special{pa 4670 1850}%
\special{pa 4650 1917}%
\special{pa 4670 1903}%
\special{pa 4690 1917}%
\special{pa 4670 1850}%
\special{fp}%
%
\special{pn 20}%
\special{pa 3400 1770}%
\special{pa 3410 1770}%
\special{fp}%
\special{sh 1}%
\special{pa 3410 1770}%
\special{pa 3343 1750}%
\special{pa 3357 1770}%
\special{pa 3343 1790}%
\special{pa 3410 1770}%
\special{fp}%
\special{pa 4850 1770}%
\special{pa 4860 1770}%
\special{fp}%
\special{sh 1}%
\special{pa 4860 1770}%
\special{pa 4793 1750}%
\special{pa 4807 1770}%
\special{pa 4793 1790}%
\special{pa 4860 1770}%
\special{fp}%
\special{pa 1950 1770}%
\special{pa 1960 1770}%
\special{fp}%
\special{sh 1}%
\special{pa 1960 1770}%
\special{pa 1893 1750}%
\special{pa 1907 1770}%
\special{pa 1893 1790}%
\special{pa 1960 1770}%
\special{fp}%
\put(48.5000,-9.6000){\makebox(0,0)[lb]{$x_1$}}%
\put(34.0000,-9.6000){\makebox(0,0)[lb]{$x_j$}}%
\put(19.5000,-9.6000){\makebox(0,0)[lb]{$x_g$}}%
\put(51.0000,-20.0000){\makebox(0,0)[lb]{$y_1$}}%
\put(36.5000,-20.0000){\makebox(0,0)[lb]{$y_j$}}%
\put(22.0000,-20.0000){\makebox(0,0)[lb]{$y_g$}}%
\end{picture}%
\caption{}
\label{figure}
\end{figure}
\end{center}
\end{tadothm}

By using the harmonic volume of the compact Riemann surface $X$ 
whose coefficients are extended over $\C$, 
Harris \cite{H-3} studied the problem of characterizing the condition 
when the cycles $W_1$ and $W_1^-$ are algebraically equivalent to each other.
Here $W_1$ is the image of the Abel-Jacobi map $X\to J(X)$ and $W_1^-$ is the image of $W_1$ under the involution $(-1)$ of $J(X)$.
Harmonic volumes or extended ones tell us the non-triviality of $W_1-W_1^-$ in $J(X)$ as follows.
If $W_1-W_1^-$ is trivial as an algebraic cycle,
 then $2I\equiv 0$ modulo $\Z$.
As is well known, if $X$ is hyperelliptic, 
then $W_1-W_1^-$ is trivial.
It is known that $I \equiv 0$ or $I \equiv 1/{2}$ modulo $\Z$ for any hyperelliptic curve $C$ by the hyperelliptic involution. 
It has been still unknown which elements in $(H^{\otimes 3})^\prime$
have nontrivial $I$ or not.
Our main theorem gives the complete answer for
this problem.

We have two ways to compute
 the harmonic volumes of all the hyperelliptic curves
in Theorem \ref{main theorem}.
One is an analytic way and the other is a topological.
In the first way, the computation of the harmonic volumes of all
the hyperelliptic curves can be reduced to that of a single hyperelliptic
 curve $C_0$, which is considered as a point of
the moduli space of hyperelliptic curves, denoted by
 $\mathcal{H}_g$.
The harmonic volume $I$ varies continuously
on the whole Torelli space $\mathcal{I}_g$,
which is the space consisting of all the compact Riemann surfaces
 with a fixed symplectic basis of $H$.
Gunning \cite{G} obtained quadratic periods
of hyperelliptic curves.
The periods are defined by iterated integrals of
 holomorphic $1$-forms along loops.
In general, iterated integrals are not homotopy invariant
 with fixed endpoints.
When we add some correction terms, they are homotopy invariant.
Because of the correction terms, the computation of harmonic volumes
 is more difficult than that of quadratic periods.
In the second way, we use basic results on the cohomology of
the hyperelliptic mapping class group.
It is denoted by $\Delta_g$.
The following theorem is obtained in the second topological way.
\begin{tadotadothm}
We have
$$\Hom_{\Delta_g}((H^{\otimes 3})^\prime, \Z/{2\Z})\cong \Z/{2\Z}.$$
\end{tadotadothm}

This theorem gives a geometric interpretation of a theorem
 established by Tanaka \cite{T}.
 It is concerning about the first homology group of $\Delta_g$
 with coefficients in $H$.

\begin{tanakathm}
$(\mathrm{Tanaka}$\cite{T}, $\mathrm{Theorem}$ $1.1)$\\
If $g\geq 2$, then
$$H_1(\Delta_g; H)=\Z/{2\Z}.$$
\end{tanakathm}

We denote by $\delta$ a connected homomorphism $H^0(\Delta_g; ((H^{\otimes 3})^\prime)^\ast)\to H^1(\Delta_g; H^\ast)^{\oplus 3}$ defined in Section $5$.
We may regard the restriction of $\delta I|_H$ as the generator of $H_1(\Delta_g; H)$.

\noindent
{\bf Acknowledgments.}
The author is grateful to Nariya Kawazumi
for valuable advice and support,
who has been encouraging the author for a long time.
He would like to thank Takeshi Katsura and
Tetsushi Ito for their helpful discussions and useful comments.
He also would like to thank the referee for useful comments.

\tableofcontents

\section{Preliminaries}
 In this section, we define a harmonic volume of a compact Riemann surface $X$ of genus $g \geq 3$. 
We begin with recalling the definition of an iterated integral on $X$. Let $\gamma: [0,1] \to X$ be a path in $X$, and $A^1(X)$ the $1$-forms on $X$. 
The iterated integral of $1$-forms $\omega_1, \omega_2, \ldots, \omega_k \in A^1(X)$ along  $\gamma$ is defined by
$$\int_{\gamma}\omega_1\omega_2\cdots\omega_k =\int_{0\leq t_1\leq t_2\leq \cdots \leq t_k \leq 1}f_1(t_1)f_2(t_2)\cdots f_k(t_k)dt_1 dt_2 \cdots dt_k,$$
where $\gamma^\ast(\omega_i)=f_i(t)dt$ in terms of the coordinate $t$ on the interval $[0,1]$. The integral is not invariant under homotopy 
with fixed endpoints.
 But, the following lemma is well known. See Chen \cite{C} for details.

\begin{lem}\label{fundamental lemma}\ \\
Let $\omega_{1, i}, \omega_{2, i}, i=1,2,\ldots, m$ be closed $1$-forms on $X$ and $\gamma$ a path in $X$. Suppose that $\displaystyle \sum_{i=1}^m \int_X \omega_{1, i} \wedge \omega_{2, i}=0$.
Take a $1$-form $\eta$ on $X$ satisfying $\displaystyle d\eta =\sum_{i=1}^m\omega_{1, i} \wedge \omega_{2, i}$. Then the integral
$$ \sum_{i=1}^m\int_{\gamma}\omega_{1, i}\omega_{2, i}-\int_{\gamma} \eta$$
is invariant under homotopy with fixed endpoints.
\end{lem}

Using iterated integrals, Harris \cite{H-1} defined the harmonic volume in the following way. In order to define it, we have to define a pointed harmonic volume for $(X, x_0)$, where $x_0$ is a point on $X$.
 We identify $H_1(X;\mathbb{Z})$ with $H^1(X;\mathbb{Z})$ by Poincar\'e duality and call them $H$. Let $K$ be the kernel of $(\ , \ ) :H\otimes H \to \Z$ induced by the intersection pairing. On the compact Riemann surface $X$, the Hodge star operator $\ast : A^1(X) \to A^1(X)$ is locally given by $\ast (f_1(z)dz + f_2(z)d\bar{z})=-\sqrt{-1}f_1(z)dz + \sqrt{-1}f_2(z)d\bar{z}$ in a local coordinate $z$ and depends only on the complex structure and not the choice of Hermitian metric. Let $\mathsf{H}_\mathbb{Z}$ denote the free abelian group of rank $2g$ spanned by all the real harmonic $1$-forms on $X$ with integral periods. We identify $H$ with $\mathsf{H}_\mathbb{Z}$ by the Hodge theorem.

\begin{defn} \mbox{(The pointed harmonic volume \cite{P})}\\
The pointed harmonic volume for $(X, x_0)$ is a linear form on $K \otimes H$ with values in $\R/{\Z}$ defined by
$$I_{x_0}{\Biggl(}\sum_{k=1}^{m}{\biggl(}\sum_{i=1}^{n_k}a_{i, k}\otimes b_{i, k}{\biggr)}\otimes c_{k}{\Biggr)}=\sum_{k=1}^{m}{\biggl(}\sum_{i=1}^{n_k}\int_{\gamma_{k}}a_{i, k}b_{i, k}-\int_{\gamma_{k}}\eta_k{\biggr)} \quad \mathrm{mod} \ \mathbb{Z},$$
where $\gamma_{k}$ is a loop in $X$ with the base point $x_0$, whose homology class is Poincar\'e dual of the cohomology class of $c_k$ and $\eta_k$ is a $1$-form on $X$, which satisfies $d\eta_k =\sum_{i=1}^{n_k}a_{i, k} \wedge b_{i, k}$ and $\int_X \eta_k\wedge\ast\alpha =0$ for any closed $1$-form $\alpha$ on $X$.
\end{defn}

 The harmonic volume is given as a restriction of the pointed harmonic volume.
 A natural homomorphism $p: H^{\otimes 3} \to H^{\oplus 3}$ is defined by
 $p(a\otimes b\otimes c)=((a, b)c, (b, c)a, (c, a)b)$.
 We denote by $(H^{\otimes 3})^{\prime}$ the kernel of $p$. It is a free $\Z$ module and satisfies the following short exact sequence
$$\xymatrix{
0 \ar[r]& (H^{\otimes 3})^{\prime} \ar[r]& H^{\otimes 3} \ar[r]^{p} & H^{\oplus 3} \ar[r] & 0 
}.$$
The rank of $(H^{\otimes 3})^{\prime}$ is $(2g)^3-6g$ and
$(H^{\otimes 3})^{\prime} \subset K \otimes H$. Harris \cite{H-1} proved that the restriction of the pointed harmonic volume on $K\otimes H$ to $(H^{\otimes 3})^{\prime}$ is independent of the choice of the base point.
\begin{defn} \mbox{(The harmonic volume \cite{H-1})}\\
The harmonic volume $I$ for $X$ is a linear form on $(H^{\otimes 3})^{\prime}$ with values in $\R/{\Z}$ defined by
$$I{\Big(}\sum_i a_i\otimes b_i\otimes c_i{\Big)}=I_{x_0}{\Big(}\sum_i a_i\otimes b_i\otimes c_i{\Big)}\quad \mathrm{mod}\ \Z.$$
\end{defn}

The map $I$ is a well-defined homomorphism 
$(H^{\otimes 3})^{\prime} \to \R/{\Z}$.
We have 
$I(\sum_i h_{\sigma(1), i}\otimes h_{\sigma(2), i}\otimes h_{\sigma(3), i})
=\mathrm{sgn}(\sigma)I(\sum_i h_{1, i}\otimes h_{2, i}\otimes h_{3, i})$, 
where $\sum_i h_{1, i}\otimes h_{2, i}\otimes h_{3, i}
 \in (H^{\otimes 3})^{\prime}$
 and $\sigma$ is an element of the third symmetric group $S_3$.
See Harris (Lemma 2.7 in \cite{H-1}) and Pulte \cite{P} for details.
In the sequel, we regard $(H^{\otimes 3})^{\prime}$ as an $S_3$-module
by this action.
 We choose a symplectic basis $\{x_i, y_i\}_{i=1, 2,  \ldots, g}$ of $H$ so that $(x_i, x_j)=(y_i, y_j)=0$ and
$(x_i, y_j)=\delta_{ij}=-(y_j, x_i)$, where $\delta_{ij}$ is Kronecker's delta. Let $z_i$ denote $x_i$ or $y_i$. We define the subset $\mathfrak{A} \subset (H^{\otimes 3})^{\prime}$ consisting of the following elements,

\begin{tabular}{lll}
$(1$ $)$ & $z_i\otimes z_j\otimes z_k$ & \hspace{-50pt}$(i\neq j\neq k\neq i)$\\
$(2$ $)$ & $x_i\otimes y_i\otimes z_k -x_{k+1}\otimes y_{k+1}\otimes z_k$ & \hspace{-50pt}$(i\neq k$ and $i\neq k+1)$ \\
$(3$a$)$ & $x_i\otimes x_i\otimes z_k $ & \hspace{-50pt}$(i\neq k)$\\
$(3$b$)$ & $y_i\otimes y_i\otimes z_k$ & \hspace{-50pt}$(i\neq k)$\\
$(4$a$)$ & $x_i\otimes x_i\otimes x_i$ & \\
$(4$b$)$ & $y_i\otimes y_i\otimes y_i$ & \\
$(5$a$)$ & $x_{i+1}\otimes x_i\otimes y_{i+1} +y_{i+1}\otimes x_i\otimes x_{i+1}$\\
$(5$b$)$ & $y_{i+1}\otimes y_i\otimes x_{i+1} +x_{i+1}\otimes y_i\otimes y_{i+1}$ \\
$(6$a$)$ & $x_i\otimes x_i\otimes y_i -x_i\otimes x_{i+1}\otimes y_{i+1} -x_{i+1}\otimes x_i\otimes y_{i+1}$ &\\
$(6$b$)$ & $y_i\otimes y_i\otimes x_i -y_i\otimes y_{i+1}\otimes x_{i+1} -y_{i+1}\otimes y_i\otimes x_{i+1}.$ & \\
\end{tabular}\\
 Here $i,j,k \in \{1,2,\ldots, g\}$ and all subscripts are read modulo $g$.
Then $\mathfrak{B}=\{\sigma (a) ; a\in \mathfrak{A}, \sigma\in S_3\}$ is a basis of $(H^{\otimes 3})^{\prime}$.

By the definition of the harmonic volume, it is obvious that $I  = 0\ \mathrm{mod}\ \Z$ for the type $(3), (4)$ and $(5)$. Furthermore, $I = 1/{2}\ \mathrm{mod}\ \Z$ for the type $(6)$. So it is enough to consider the type $(1)$ and $(2)$.

\section{The periods and iterated integrals of a hyperelliptic curve}
In this section, we compute the periods and iterated integrals of a hyperelliptic curve of genus $g\geq 3$. First of all, we take a symplectic basis of $H$.
\subsection{A homology basis of hyperelliptic curves}\ \\
We define a hyperelliptic curve $C$ as follows. Let $p_0, p_1,\ldots, p_{2g+1}$ be distinct points on $\C$. It is the compactification of the plane curve in the $(z,w)$ plane $\C^2$
$$w^2=\prod_{i=0}^{2g+1} (z-p_i),$$
and admits the hyperelliptic involution given by $\iota:(z,w) \mapsto (z,-w)$.
Let $\pi$ be the $2$-sheeted covering $C \to \C P^1, (z,w) \mapsto z$, branched over $2g+2$ branch points $\{p_i\}_{i=0,1,\cdots ,2g+1}$ and $P_i \in C$
a ramification point so that $\pi(P_i)=p_i$.
 On the curve $C$, we choose endpoints $Q_0, Q_1(=\iota(Q_0))$
 as in Figure \ref{path}.
We define by $\Omega$ the simply-connected domain
 $\displaystyle\C P^1\setminus \bigcup_{j=0}^{g}p_{2j}p_{2j+1}$,
where $p_{2j}p_{2j+1}$ is a simple arc
connecting $p_{2j}$ and $p_{2j+1}$.
 Then $\pi^{-1}(\Omega)$ consists of two connected components.
 We denote by $\Omega_0,\Omega_1$ the connected components of
 $\pi^{-1}(\Omega)$ which contain $Q_0,Q_1$ respectively. 
Let $e_j, j=0,1, \ldots, 2g+1$, be a path in $C$ which is to be followed from $Q_0$ to $P_j$ and go to $Q_1$ along the arcs $Q_0P_j$ and $P_jQ_1$.
 See Figure \ref{path}.
We write simply $\overline{e}_j$ for $\pi(e_j)$.
 It is a loop in $\C P^1$ with the base point $\pi(Q_0)$.
\begin{center}
\begin{figure}[htbp]
\vspace*{10pt}

\unitlength 0.1in
\begin{picture}(56.00,16.10)(6.80,-22.00)
%
\special{pn 13}%
\special{pa 3400 2040}%
\special{pa 3419 2014}%
\special{pa 3438 1988}%
\special{pa 3456 1962}%
\special{pa 3475 1935}%
\special{pa 3493 1909}%
\special{pa 3510 1882}%
\special{pa 3528 1855}%
\special{pa 3544 1828}%
\special{pa 3560 1800}%
\special{pa 3575 1772}%
\special{pa 3589 1744}%
\special{pa 3603 1715}%
\special{pa 3616 1686}%
\special{pa 3628 1657}%
\special{pa 3640 1627}%
\special{pa 3652 1597}%
\special{pa 3663 1567}%
\special{pa 3674 1537}%
\special{pa 3684 1506}%
\special{pa 3695 1476}%
\special{pa 3705 1445}%
\special{pa 3715 1415}%
\special{pa 3720 1400}%
\special{sp}%
%
\special{pn 20}%
\special{pa 1800 600}%
\special{pa 5000 600}%
\special{fp}%
\special{pa 5000 2200}%
\special{pa 1800 2200}%
\special{fp}%
%
\special{pn 20}%
\special{ar 1800 1400 800 800  1.5707963 4.7123890}%
%
\special{pn 20}%
\special{ar 5000 1400 800 800  4.7123890 6.2831853}%
\special{ar 5000 1400 800 800  0.0000000 1.5707963}%
%
\special{pn 20}%
\special{ar 1960 1360 280 160  6.2831853 6.2831853}%
\special{ar 1960 1360 280 160  0.0000000 3.1415927}%
%
\special{pn 20}%
\special{ar 1960 1400 248 104  3.1415927 6.2831853}%
%
\special{pn 20}%
\special{ar 3400 1360 280 160  6.2831853 6.2831853}%
\special{ar 3400 1360 280 160  0.0000000 3.1415927}%
%
\special{pn 20}%
\special{ar 4840 1360 280 160  6.2831853 6.2831853}%
\special{ar 4840 1360 280 160  0.0000000 3.1415927}%
%
\special{pn 20}%
\special{ar 3400 1400 248 104  3.1415927 6.2831853}%
%
\special{pn 20}%
\special{ar 4840 1400 248 104  3.1415927 6.2831853}%
%
\special{pn 8}%
\special{sh 0.600}%
\special{ar 3400 2040 24 24  0.0000000 6.2831853}%
%
\special{pn 8}%
\special{sh 0.600}%
\special{ar 1040 1400 24 24  0.0000000 6.2831853}%
%
\special{pn 8}%
\special{sh 0.600}%
\special{ar 5760 1400 24 24  0.0000000 6.2831853}%
%
\special{pn 8}%
\special{sh 0.600}%
\special{ar 5144 1400 24 24  0.0000000 6.2831853}%
%
\special{pn 8}%
\special{sh 0.600}%
\special{ar 3704 1400 24 24  0.0000000 6.2831853}%
%
\special{pn 8}%
\special{sh 0.600}%
\special{ar 2264 1400 24 24  0.0000000 6.2831853}%
%
\special{pn 8}%
\special{sh 0.600}%
\special{ar 1656 1400 24 24  0.0000000 6.2831853}%
%
\special{pn 8}%
\special{sh 0.600}%
\special{ar 3096 1400 24 24  0.0000000 6.2831853}%
%
\special{pn 8}%
\special{sh 0.600}%
\special{ar 4536 1400 24 24  0.0000000 6.2831853}%
%
\special{pn 8}%
\special{sh 0.600}%
\special{ar 3400 760 24 24  0.0000000 6.2831853}%
%
\special{pn 20}%
\special{sh 1}%
\special{ar 2680 1400 10 10 0  6.28318530717959E+0000}%
\special{sh 1}%
\special{ar 2600 1400 10 10 0  6.28318530717959E+0000}%
\special{sh 1}%
\special{ar 2760 1400 10 10 0  6.28318530717959E+0000}%
%
\special{pn 20}%
\special{sh 1}%
\special{ar 4120 1400 10 10 0  6.28318530717959E+0000}%
\special{sh 1}%
\special{ar 4040 1400 10 10 0  6.28318530717959E+0000}%
\special{sh 1}%
\special{ar 4200 1400 10 10 0  6.28318530717959E+0000}%
%
\special{pn 8}%
\special{pa 1000 1400}%
\special{pa 680 1400}%
\special{fp}%
\special{pa 5800 1400}%
\special{pa 6280 1400}%
\special{fp}%
%
\special{pn 13}%
\special{pa 3720 1400}%
\special{pa 3710 1369}%
\special{pa 3700 1339}%
\special{pa 3689 1308}%
\special{pa 3679 1278}%
\special{pa 3668 1247}%
\special{pa 3657 1217}%
\special{pa 3646 1187}%
\special{pa 3634 1158}%
\special{pa 3622 1128}%
\special{pa 3609 1099}%
\special{pa 3596 1070}%
\special{pa 3582 1041}%
\special{pa 3567 1013}%
\special{pa 3552 985}%
\special{pa 3535 958}%
\special{pa 3519 931}%
\special{pa 3501 904}%
\special{pa 3483 877}%
\special{pa 3465 851}%
\special{pa 3447 825}%
\special{pa 3428 798}%
\special{pa 3409 772}%
\special{pa 3400 760}%
\special{sp -0.045}%
%
\special{pn 20}%
\special{pa 3600 1720}%
\special{pa 3608 1704}%
\special{fp}%
\special{sh 1}%
\special{pa 3608 1704}%
\special{pa 3560 1755}%
\special{pa 3584 1752}%
\special{pa 3596 1773}%
\special{pa 3608 1704}%
\special{fp}%
%
\special{pn 20}%
\special{pa 3600 1080}%
\special{pa 3592 1064}%
\special{fp}%
\special{sh 1}%
\special{pa 3592 1064}%
\special{pa 3604 1133}%
\special{pa 3616 1112}%
\special{pa 3640 1115}%
\special{pa 3592 1064}%
\special{fp}%
%
\special{pn 8}%
\special{ar 6040 1400 160 480  3.7310362 6.2831853}%
\special{ar 6040 1400 160 480  0.0000000 2.3561945}%
%
\special{pn 8}%
\special{pa 5928 1720}%
\special{pa 5912 1680}%
\special{fp}%
\special{sh 1}%
\special{pa 5912 1680}%
\special{pa 5918 1749}%
\special{pa 5932 1730}%
\special{pa 5955 1734}%
\special{pa 5912 1680}%
\special{fp}%
\put(34.8000,-7.6000){\makebox(0,0)[lb]{$Q_1$}}%
\put(34.8000,-21.2000){\makebox(0,0)[lb]{$Q_0$}}%
\put(56.0000,-16.0000){\makebox(0,0)[lb]{$P_0$}}%
\put(11.2000,-16.0000){\makebox(0,0)[lb]{$P_{2g+1}$}}%
\put(14.8000,-16.0000){\makebox(0,0)[lb]{$P_{2g}$}}%
\put(29.2000,-16.0000){\makebox(0,0)[lb]{$P_{j+1}$}}%
\put(43.6000,-16.0000){\makebox(0,0)[lb]{$P_2$}}%
\put(52.0000,-16.0000){\makebox(0,0)[lb]{$P_1$}}%
\put(37.6000,-16.0000){\makebox(0,0)[lb]{$P_j$}}%
\put(23.2000,-16.0000){\makebox(0,0)[lb]{$P_{2g-1}$}}%
\put(58.4000,-8.8000){\makebox(0,0)[lb]{$180^\circ$}}%
\put(60.4000,-20.4000){\makebox(0,0)[lb]{$\iota$}}%
\put(58.0000,-22.0000){\makebox(0,0)[lb]{$C$}}%
\end{picture}%

\vspace*{15pt}

{\Large $\downarrow^\pi$}

\vspace*{15pt}

\hspace*{-41pt}
\unitlength 0.1in
\begin{picture}( 51.5000, 11.2000)(  6.5000,-26.0000)
%
\special{pn 13}%
\special{pa 3400 2440}%
\special{pa 3420 2414}%
\special{pa 3438 2388}%
\special{pa 3456 2362}%
\special{pa 3476 2336}%
\special{pa 3494 2310}%
\special{pa 3510 2282}%
\special{pa 3528 2256}%
\special{pa 3544 2228}%
\special{pa 3560 2200}%
\special{pa 3576 2172}%
\special{pa 3590 2144}%
\special{pa 3604 2116}%
\special{pa 3616 2086}%
\special{pa 3628 2058}%
\special{pa 3640 2028}%
\special{pa 3652 1998}%
\special{pa 3664 1968}%
\special{pa 3674 1938}%
\special{pa 3684 1906}%
\special{pa 3696 1876}%
\special{pa 3706 1846}%
\special{pa 3716 1816}%
\special{pa 3720 1800}%
\special{sp}%
%
\special{pn 20}%
\special{ar 1800 1800 800 800  1.5707963 3.1415927}%
%
\special{pn 20}%
\special{ar 5000 1800 800 800  6.2831853 6.2831853}%
\special{ar 5000 1800 800 800  0.0000000 1.5707963}%
%
\special{pn 20}%
\special{ar 1960 1760 280 160  1.5707963 2.7851557}%
%
\special{pn 20}%
\special{ar 3400 1760 280 160  0.3290497 2.8399736}%
%
\special{pn 20}%
\special{ar 4840 1760 280 160  0.3472408 1.5707963}%
%
\special{pn 8}%
\special{sh 0.600}%
\special{ar 3400 2440 24 24  0.0000000 6.2831853}%
%
\special{pn 8}%
\special{sh 0.600}%
\special{ar 1040 1800 24 24  0.0000000 6.2831853}%
%
\special{pn 8}%
\special{sh 0.600}%
\special{ar 5760 1800 24 24  0.0000000 6.2831853}%
%
\special{pn 8}%
\special{sh 0.600}%
\special{ar 5144 1800 24 24  0.0000000 6.2831853}%
%
\special{pn 8}%
\special{sh 0.600}%
\special{ar 3704 1800 24 24  0.0000000 6.2831853}%
%
\special{pn 8}%
\special{sh 0.600}%
\special{ar 1656 1800 24 24  0.0000000 6.2831853}%
%
\special{pn 8}%
\special{sh 0.600}%
\special{ar 3096 1800 24 24  0.0000000 6.2831853}%
%
\special{pn 20}%
\special{sh 1}%
\special{ar 2400 1920 10 10 0  6.28318530717959E+0000}%
\special{sh 1}%
\special{ar 2320 1920 10 10 0  6.28318530717959E+0000}%
\special{sh 1}%
\special{ar 2480 1920 10 10 0  6.28318530717959E+0000}%
%
\special{pn 20}%
\special{sh 1}%
\special{ar 4400 1920 10 10 0  6.28318530717959E+0000}%
\special{sh 1}%
\special{ar 4320 1920 10 10 0  6.28318530717959E+0000}%
\special{sh 1}%
\special{ar 4480 1920 10 10 0  6.28318530717959E+0000}%
%
\special{pn 20}%
\special{pa 3600 2120}%
\special{pa 3608 2104}%
\special{fp}%
\special{sh 1}%
\special{pa 3608 2104}%
\special{pa 3560 2156}%
\special{pa 3584 2152}%
\special{pa 3596 2174}%
\special{pa 3608 2104}%
\special{fp}%
\put(34.8000,-25.2000){\makebox(0,0)[lb]{$\pi(Q_0)=\pi(Q_1)$}}%
\put(56.0000,-16.5000){\makebox(0,0)[lb]{$p_0$}}%
\put(11.2000,-16.5000){\makebox(0,0)[lb]{$p_{2g+1}$}}%
\put(14.8000,-16.5000){\makebox(0,0)[lb]{$p_{2g}$}}%
\put(29.2000,-16.5000){\makebox(0,0)[lb]{$p_{j+1}$}}%
\put(52.0000,-16.5000){\makebox(0,0)[lb]{$p_1$}}%
\put(37.6000,-16.5000){\makebox(0,0)[lb]{$p_j$}}%
\put(58.0000,-26.0000){\makebox(0,0)[lb]{$\mathbb{C} P^1$}}%
%
\special{pn 20}%
\special{pa 1800 2600}%
\special{pa 5000 2600}%
\special{fp}%
%
\special{pn 20}%
\special{ar 1350 1800 280 110  0.0000000 6.2831853}%
%
\special{pn 20}%
\special{ar 5450 1800 280 110  0.0000000 6.2831853}%
%
\special{pn 20}%
\special{ar 2780 1800 280 110  4.7123890 6.2831853}%
\special{ar 2780 1800 280 110  0.0000000 1.5707963}%
%
\special{pn 20}%
\special{ar 4020 1800 280 110  1.5707963 4.7123890}%
%
\special{pn 13}%
\special{pa 3700 1800}%
\special{pa 3680 1826}%
\special{pa 3660 1852}%
\special{pa 3642 1878}%
\special{pa 3622 1904}%
\special{pa 3604 1930}%
\special{pa 3584 1956}%
\special{pa 3566 1982}%
\special{pa 3550 2010}%
\special{pa 3534 2036}%
\special{pa 3518 2064}%
\special{pa 3504 2092}%
\special{pa 3490 2122}%
\special{pa 3478 2150}%
\special{pa 3468 2180}%
\special{pa 3458 2210}%
\special{pa 3448 2242}%
\special{pa 3440 2272}%
\special{pa 3430 2304}%
\special{pa 3422 2334}%
\special{pa 3416 2366}%
\special{pa 3408 2398}%
\special{pa 3400 2430}%
\special{pa 3400 2430}%
\special{sp}%
%
\special{pn 20}%
\special{pa 3490 2120}%
\special{pa 3480 2140}%
\special{fp}%
\special{sh 1}%
\special{pa 3480 2140}%
\special{pa 3528 2090}%
\special{pa 3504 2092}%
\special{pa 3492 2072}%
\special{pa 3480 2140}%
\special{fp}%
\put(37.0000,-21.7000){\makebox(0,0)[lb]{$\overline{e}_j=\pi(e_j)$}}%
%
\special{pn 8}%
\special{pa 1000 1800}%
\special{pa 650 1800}%
\special{ip}%
\end{picture}%

\vspace*{20pt}
\caption{}
\label{path}
\end{figure}
\end{center}

It is obvious that $e_{j_1}\cdot \iota(e_{j_2})$ is a loop in $C$ with the base point $Q_0$, where the product $e_{j_1}\cdot \iota(e_{j_2})$ indicates that we traverse $e_{j_1}$ first, then $\iota(e_{j_2})$.
So we have the homotopy equivalences relative to the base point $Q_0$
$$e_j\cdot \iota(e_j) \sim 1, \quad j=0,1, \ldots ,2g+1,$$
and
$$e_0 \cdot \iota(e_1)\cdot \cdots \cdot e_{2g} \cdot \iota(e_{2g+1}) \sim 1.$$
We define the loops $a_i, b_i, i=1,2,\ldots, g$, in $C$
 with the base point $Q_0$ by 
$$
\left.\begin{array}{l}
\displaystyle a_i=e_{2i-1}\cdot\iota(e_{2i}) ,\\
\displaystyle b_i=e_{2i-1}\cdot\iota(e_{2i-2})\cdot\cdots\cdot e_1\cdot\iota(e_0).
\end{array}\right.
$$
So a symplectic basis of $H_1(C; \Z)$ can be given by 
$\{[a_i], [b_i]\}_{i=1,2,\ldots, g}$, where $[a_i]$ and $[b_i]$
 are the homology classes of $a_i$ and $b_i$ respectively.
In fact, we have $([a_i], [b_j])=\delta_{ij}=-([b_j], [a_i])$ and 
$([a_i], [a_j])=([b_i], [b_j])=0$.
It is clear that $[a_i]$ and $[b_i]$ are
equal to $x_i$ and $y_i$ in Figure \ref{figure} respectively.

\subsection{The hyperelliptic curve $C_0$}\ \\
A hyperelliptic curve $C_0$ is defined by the equation $w^2=z^{2g+2}-1$. We take $Q_i=(0, (-1)^{i}\sqrt{-1}), i=0, 1$, and $P_j=(\zeta^j, 0), j=0,1, \ldots, 2g+1$, where $\zeta=e^{2\pi\sqrt{-1}/(2g+2)}$. We define a path $e_j :[0,1] \to C_0, j=0, 1,\ldots ,2g+1$, by 
$$
\left\{
  \begin{array}{lll}
  \big( 2t\zeta^j,\sqrt{-1}\sqrt{1-(2t)^{2g+2}} \big)& \mathrm{for} &0\leq t\leq 1/{2}, \\
  \big( (2-2t)\zeta^j,-\sqrt{-1}\sqrt{1-(2-2t)^{2g+2}} \big)& \mathrm{for} & 1/{2}\leq t\leq 1.
  \end{array}
  \right.
$$

We denote $\omega_i=z^{i-1}dz/{w}, i=1,2, \ldots, g$,
which are holomorphic $1$-forms on $C_0$. 
Then $\{\omega_i\}_{i=1,2,\ldots ,g}$ is a basis of the space of holomorphic
 $1$-forms on $C_0$.
 Let $B(u,v)$ denote the beta function 
$\int_0^1x^{u-1}(1-x)^{v-1}dx$ for $u,v>0$. It is easy to show.
\begin{lem}
\label{lem1}
We have
$$\int_{e_j}\omega_i=-2\sqrt{-1}\zeta^{ij}B(i/(2g+2),1/2)/{(2g+2)}=-\int_{\iota(e_j)}\omega_i.$$
\end{lem}

We denote by $\omega^\prime_i$ the holomorphic $1$-form $\displaystyle {(2g+2)\sqrt{-1}\over{2B(i/(2g+2),1/2)}}\omega_i$. The periods of $C_0$ are obtained by Lemma \ref{lem1}.
\begin{lem}\label{lem2}
We have
$$
\left.\begin{array}{l}
\displaystyle \int_{a_j}\omega^\prime_i=\zeta^{i(2j-1)}(1-\zeta^i),\\
\displaystyle \int_{b_j}\omega^\prime_i={{\zeta^{2ij}-1}\over{\zeta^i}+1},
\end{array}\right.$$
where $i, j \in \{1,2, \ldots, g\}$.
\end{lem}
\begin{rem}
Since $\omega^\prime_i$ is a closed $1$-form, the integral $\int_{\gamma} \omega^\prime_i$ depends only on the homology classes of $\gamma$.
\end{rem}

In order to prove Lemma \ref{iterated integrals},
we start with the following well known lemma.
\begin{lem}\label{shuffle product}
Let $\omega_1,\omega_2$ be $1$-forms on $X$
and $\gamma_1,\gamma_2,\ldots,\gamma_m$ paths in $X$
so that $\gamma_1\gamma_2\cdots \gamma_m$ is a path.
Then, we have
$$
\int_{\gamma_1\gamma_2\cdots \gamma_m}\omega_1\omega_2
=\sum_{i=1}^{m}\int_{\gamma_i}\omega_1\omega_2
+\sum_{i<j}\int_{\gamma_i}
\omega_1\int_{\gamma_j}\omega_2.
$$
\end{lem}
Since $\iota$ is a diffeomorphism of $C_0$ and $\iota(e_k)=e_k^{-1}$,
 we have
$$\int_{e_k}\omega^{\prime}_i\omega^{\prime}_j
=\int_{\iota(e_k)}\omega^{\prime}_i\omega^{\prime}_j
=\int_{e_k^{-1}}\omega^{\prime}_i\omega^{\prime}_j
=-\int_{e_k}\omega^{\prime}_i\omega^{\prime}_j-
\int_{e_k}\omega^{\prime}_i\int_{e_k^{-1}}\omega^{\prime}_j.$$
Then $\displaystyle
\int_{e_k}\omega^{\prime}_i\omega^{\prime}_j=
\frac{1}{2}\int_{e_k}\omega^{\prime}_i\int_{e_k}\omega^{\prime}_j$.
This formula, Lemma \ref{lem2} and Lemma \ref{shuffle product}
give us iterated integrals of $\omega^\prime_i$ along $a_k$ and $b_k$.

\begin{lem}\label{iterated integrals}
We have
$$
\left.\begin{array}{l}
\displaystyle \int_{a_k}\omega^\prime_i\omega^\prime_j={1\over2}\zeta^{(i+j)(2k-1)}(1-2\zeta^j+\zeta^{i+j}),\\
\displaystyle \int_{b_k}\omega^\prime_i\omega^\prime_j=\sum_{l=1}^k{1\over2}\zeta^{(i+j)(2l-2)}(1-2\zeta^i+\zeta^{i+j})+\sum_{1\leq l<m\leq k}(\zeta^i-1)(\zeta^j-1)\zeta^{i(2m-2)+j(2l-2)},
\end{array}\right.
$$
where $i, j \in \{1,2, \ldots, g\}$.
\end{lem}

For the rest of this section, we compute the iterated integrals of real harmonic $1$-forms  of $C_0$ with integral periods. Let $\Omega_a$ and $\Omega_b$ be the non-singular matrices \\
$$\left(
\begin{array}{ccc}
\int_{a_1}\omega_1^\prime & \cdots & \int_{a_g}\omega_1^\prime \\
\vdots & & \vdots \\
\int_{a_1}\omega_g^\prime & \cdots & \int_{a_g}\omega_g^\prime
\end{array}
\right)\ \mathrm{and}\ 
\left(
\begin{array}{ccc}
\int_{b_1}\omega_1^\prime & \cdots & \int_{b_g}\omega_1^\prime \\
\vdots & & \vdots \\
\int_{b_1}\omega_g^\prime & \cdots & \int_{b_g}\omega_g^\prime
\end{array}
\right),$$
respectively. 
It is clear that $(ij)$-entries of $(\Omega_a)^{-1}$ and $(\Omega_b)^{-1}$
are given by 
$\displaystyle\frac{1}{g+1}
\frac{\zeta^j(-1+\zeta^{-2ij})}{1-\zeta^j}$ and 
$\displaystyle\frac{1}{g+1}
\zeta^{-2ij}(1+\zeta^j)
$ respectively. 
Then we obtain the period matrix $(\Omega_a)^{-1}\Omega_b$ denoted by $Z$. In general, it is well known that $Z \in GL(g, \C)$ is symmetric and its imaginary part $\Im Z$ is positive definite. In particular, Schindler \cite{S} proved the theorem below. We deduce it directly from Lemma \ref{lem2}.
\begin{thm}\label{period matrix}
$(\mathrm{Schindler}$\cite{S}$,\ \mathrm{Theorem}\ 2)$\\
Let $Z$ be the period matrix on the curve $C_0$ as above. Then its $(ij)$-entry is given by 
$$\displaystyle {1\over{g+1}}\sum_{k=1}^{g}{{\zeta^k(\zeta^{-2ik}-1)(\zeta^{2kj}-1)}\over{1-\zeta^{2k}}}.$$
Furthermore, all the entries are pure imaginary.
\end{thm}
\begin{rem}
We need some steps for another presentation of $Z$
by Schindler as follows.
\begin{align*}
\sum_{k=1}^{g}{{\zeta^k(\zeta^{-2ik}-1)(\zeta^{2kj}-1)}\over{1-\zeta^{2k}}}
&=\sum_{k=1}^{g}\zeta^k(\zeta^{2kj}-1)
\zeta^{-2k}\frac{1-\big(\zeta^{-2k}\big)^i}{1-\zeta^{-2k}}\\
&=\sum_{k=1}^{g}\zeta^k(\zeta^{2kj}-1)\sum_{\nu=1}^{i}\zeta^{-2k\nu}\\
&=\sum_{\nu=1}^{i}\sum_{k=1}^{g}\Big(\big(\zeta^{1-2\nu+2j}\big)^k
-\big(\zeta^{1-2\nu}\big)^k\Big)\\
&=\sum_{\nu=1}^{i}\biggl(
\frac{1+\zeta^{1-2\nu+2j}}{1-\zeta^{1-2\nu+2j}}
-\frac{1+\zeta^{1-2\nu}}{1-\zeta^{1-2\nu}}\biggr).
\end{align*}
Then we have
$${1\over{g+1}}\sum_{k=1}^{g}{{\zeta^k(\zeta^{-2ik}-1)
(\zeta^{2kj}-1)}\over{1-\zeta^{2k}}}
=
\frac{\sqrt{-1}}{g+1}
\Biggl(
\sum_{\nu=1}^{i}
\frac{1+\cos \frac{2\nu-1}{g+1}\pi}{\sin\frac{2\nu-1}{g+1}\pi}
+\frac{1+\cos \frac{2(j-\nu)+1}{g+1}\pi}{\sin\frac{2(j-\nu)+1}{g+1}\pi}
\Biggr).
$$
\end{rem}

We define real harmonic $1$-forms $\alpha_i, \beta_i, i=1,2, \ldots, g$, by
$$
\left(\begin{array}{c}
\alpha_1\\
 \vdots\\
\alpha_g
\end{array}
\right)=
(\Im Z)^{-1}\Im
\left(\begin{array}{c}
(\Omega_a)^{-1}
\left(\begin{array}{c}
\omega_1^\prime\\
 \vdots\\
 \omega_g^\prime
\end{array}
\right)
\end{array}
\right)
\hspace{5pt} \mathrm{and}\hspace{5pt}
\left(\begin{array}{c}
\beta_1\\
 \vdots\\
\beta_g 
\end{array}
\right)=-\Re
\left(\begin{array}{c}
(\Omega_a)^{-1}
\left(\begin{array}{c}
\omega_1^\prime\\
 \vdots\\
 \omega_g^\prime
\end{array}
\right)\end{array}
\right)
.$$
Using Theorem \ref{period matrix}, we have
$\Im Z=-\sqrt{-1}(\Omega_a)^{-1}\Omega_b$.
Then
$$
\left(\begin{array}{c}
\alpha_1\\
 \vdots\\
\alpha_g
\end{array}
\right)=
\Re
\left(\begin{array}{c}
(\Omega_b)^{-1}
\left(\begin{array}{c}
\omega_1^\prime\\
 \vdots\\
 \omega_g^\prime
\end{array}
\right)
\end{array}
\right).
$$
It is clear that $\int_{a_j}\alpha_i=\int_{b_j}\beta_i=0$ and $\int_{b_j}\alpha_i=\delta_{ij}=-\int_{a_j}\beta_i$ by Lemma \ref{lem2}. Let $\mathrm{PD}$ denote
 the Poincar\'e dual $H_1(C_0; \Z)\to H^1(C_0; \Z)$.
We have $\mathrm{PD}([a_i])=\alpha_i$ and $\mathrm{PD}([b_i])=\beta_i$ for
$i=1,2, \ldots, g$. Hence, $\{\alpha_i, \beta_i\}_{i=1,2, \ldots, g}
 \subset H^1(C_0; \Z)$ is a symplectic basis.

Let $t_u$ be a complex number
$\displaystyle \sum_{p=1}^g\zeta^{up}$ for any integer $u$.
It is obvious that
$$
t_u=
\left\{
  \begin{array}{lll}
  g & \mathrm{for} &u\in (2g+2)\Z,\\
  -1 & \mathrm{for} &u\not\in (2g+2)\Z\ \mathrm{and}\ u: \mathrm{even},\\
  \displaystyle{{1+\zeta^u}\over{1-\zeta^u}} &\mathrm{for} & u: \mathrm{odd}.
  \end{array}
  \right.
$$
Moreover, $t_u$ is pure imaginary and
$t_{-u}=-t_u$ when $u$ is odd.

Using Lemma \ref{iterated integrals}, we can calculate iterated integrals by means of $t_u$ as follows.
\begin{lem}
\label{a hyperelliptic iterated integral}
On the curve $C_0$, we have the equations
$$
\begin{array}{rcl}
(1)&\displaystyle\int_{a_k}\beta_i\beta_j &
=\displaystyle{1\over{2(g+1)^2}}{\biggl\{}(t_{2k-2j}-t_{2k})\sum_{u=1}^it_{2k-2u}+(t_{2k}-t_{2k-2i})\sum_{u=1}^jt_{2k-2u+2}{\biggr\}}, \\
(2)&\displaystyle\int_{b_k}\beta_i\beta_j &=0, \\
(3)&\displaystyle\int_{a_k}\alpha_i\alpha_j &=0, \\
(4)&\displaystyle\int_{b_k}\alpha_i\alpha_j &=
\displaystyle{1\over{2(g+1)^2}}{\biggl\{}\sum_{u=1}^k(t_{2u-2j}t_{2u-2i}-2t_{2u-2j-2}t_{2u-2i}+t_{2u-2j-2}t_{2u-2i-2})\\
&             &\hspace*{85pt} \displaystyle+\sum_{v=2}^k2(t_{2v-2i}-t_{2v-2i-2})(t_{2v-2j-2}-t_{(-2j)}){\biggr\}}.
\end{array}
$$
Here $i, j, k \in \{1,2, \ldots, g\}$.
\end{lem}
\begin{rem}
For $k=1$, $\sum_{v=2}^k2(t_{2v-2i}-t_{2v-2i-2})(t_{2v-2j-2}-t_{(-2j)})=0$.
\end{rem}
\begin{proof}
We compute $\int_{a_k}\beta_i\beta_j$ in the following way.
Let $A_{i,j}$ be the $(i,j)$-entry of $(\Omega_a)^{-1}$.
By the definition of $\beta_i$, we have
$\beta_i=-\Re\Big(\sum_{l=1}^{g}A_{i,l}\omega^{\prime}_l\Big)$.
Using this, $\int_{a_k}\beta_i\beta_j$ can be given by
\begin{align*}
&\int_{a_k}\Re\biggl
(\sum_{l=1}^{g}A_{i,l}\omega^{\prime}_l\biggr)
\Re\biggl(\sum_{m=1}^{g}A_{j,m}\omega^{\prime}_m\biggr)\\
=&
\frac{1}{4}
\int_{a_k}
\sum_{l,m=1}^{g}
\biggl(
A_{i,l}A_{j,m}\omega^{\prime}_l\omega^{\prime}_m
+A_{i,l}\overline{A}_{j,m}\omega^{\prime}_l\overline{\omega}^{\prime}_m
+\overline{A}_{i,l}A_{j,m}\overline{\omega}^{\prime}_l\omega^{\prime}_m
+\overline{A}_{i,l}\overline{A}_{j,m}
\overline{\omega}^{\prime}_l\overline{\omega}^{\prime}_m
\biggr)\\
=&
\frac{1}{2}
\Re\Biggl\{
\sum_{l,m=1}^{g}
\biggl(
A_{i,l}A_{j,m}\int_{a_k}\omega^{\prime}_l\omega^{\prime}_m
+A_{i,l}\overline{A}_{j,m}
\int_{a_k}\omega^{\prime}_l\overline{\omega}^{\prime}_m
\biggr)
\Biggr\}.
\end{align*}
Lemma \ref{iterated integrals} gives us
\begin{align*}
&(g+1)^2\sum_{l,m=1}^{g}
A_{i,l}A_{j,m}\int_{a_k}\omega^{\prime}_l\omega^{\prime}_m\\
=&
\sum_{l,m=1}^{g}
\frac{\zeta^l(-1+\zeta^{-2il})}{1-\zeta^l}
\frac{\zeta^m(-1+\zeta^{-2jm})}{1-\zeta^m}
{1\over2}\zeta^{(l+m)(2k-1)}(1-2\zeta^m+\zeta^{l+m})\\
=&\frac{1}{2}
\sum_{m=1}^{g}
\frac{1-\zeta^{2jm}}{1-\zeta^m}
\zeta^{m(2k-2j)}
\sum_{l=1}^{g}
\frac{1-\zeta^{2il}}{1-\zeta^l}
\zeta^{l(2k-2i)}(1-2\zeta^m+\zeta^{l+m})\\
=&\frac{1}{2}
\sum_{m=1}^{g}
\sum_{v=2k-2j}^{2k-1}\zeta^{mv}
\sum_{l=1}^{g}\sum_{u=2k-2i}^{2k-1}
\zeta^{lu}(1-2\zeta^m+\zeta^{l+m})\\
=&\frac{1}{2}
\sum_{m=1}^{g}
\sum_{v=2k-2j}^{2k-1}\zeta^{mv}
\biggl\{
\sum_{u=2k-2i}^{2k-1}\Big(t_u(1-2\zeta^m)+t_{u+1}\zeta^{m}\Big)
\biggr\}\\
=&\frac{1}{2}
\sum_{m=1}^{g}
\sum_{v=2k-2j}^{2k-1}\zeta^{mv}
\biggl\{
\sum_{u=2k-2i}^{2k-1}t_u(1-\zeta^m)+(t_{2k}-t_{2k-2i})\zeta^{m}
\biggr\}\\
=&\frac{1}{2}
\sum_{v=2k-2j}^{2k-1}
\sum_{m=1}^{g}
\biggl\{
\sum_{u=2k-2i}^{2k-1}
t_u\zeta^{mv}(1-\zeta^m)+(t_{2k}-t_{2k-2i})\zeta^{m(v+1)}
\biggr\}\\
=&\frac{1}{2}
\sum_{v=2k-2j}^{2k-1}
\biggl\{\sum_{u=2k-2i}^{2k-1}
t_u
(t_v-t_{v+1})+(t_{2k}-t_{2k-2i})t_{v+1}
\biggr\}.
\end{align*}
So we obtain
$$
(g+1)^2\sum_{l,m=1}^{g}
A_{i,l}A_{j,m}\int_{a_k}\omega^{\prime}_l\omega^{\prime}_m
=\frac{1}{2}
\biggl\{(t_{2k-2j}-t_{2k})\sum_{u=2k-2i}^{2k-1}
t_u+
(t_{2k}-t_{2k-2i})\sum_{v=2k-2j+1}^{2k}t_{v}
\biggr\}.$$
Using 
$\int_{a_k}\omega^\prime_i\overline{\omega}^\prime_j
={1\over2}\zeta^{(i-j)(2k-1)}(1-2\zeta^{-j}+\zeta^{i-j})$,
the value $(g+1)^2\sum_{l,m=1}^{g}
A_{i,l}\overline{A}_{j,m}\int_{a_k}\omega^{\prime}_l
\overline{\omega}^{\prime}_m$
can be computed by
$$
\frac{1}{2}
\biggl\{(t_{-(2k-2j)}-t_{-(2k)})\sum_{u=2k-2i}^{2k-1}
t_u+
(t_{2k}-t_{2k-2i})\sum_{v=2k-2j+1}^{2k}t_{-v}
\biggr\}.$$
Since $t_u=t_{-u}$ for $u\in 2\Z$ and $t_u$ is pure imaginary
 for $u\in 2\Z+1$, we have (1).
The values $\int_{b_k}\beta_i\beta_j, \int_{a_k}\alpha_i\alpha_j$
and $\int_{b_k}\alpha_i\alpha_j$
are calculated similarly.
\end{proof}

\section{The harmonic volumes of hyperelliptic curves}
In this section, we consider the harmonic volumes of hyperelliptic curve\underline{s}.
They can be reduced to the computation for the hyperelliptic curve $C_0$.

\begin{thm}\label{main theorem}
For any hyperelliptic curve $C$, let $\{x_i, y_i\}_{i=1,2, \ldots, g}$ be a symplectic basis of $H=H_1(C; \Z)$ in Figure \ref{figure},
where $\iota$ is the hyperelliptic involution.
We denote by $z_i$ either $x_i$ or $y_i$.
Then,
$$
\mbox{\phantom{gi}}I(z_i\otimes z_j \otimes z_k) =
\left.
  \begin{array}{cc}
  0 & \mbox{for}\ i\neq j\neq k\neq i,
  \end{array}
\right.
$$
$$
I(x_i\otimes y_i\otimes z_k -x_{k+1}\otimes y_{k+1}\otimes z_k)=
\left\{
  \begin{array}{cc}
  {\displaystyle {1/{2}}} & \mbox{for}\ i<k, k=2, 3, \ldots, g-1\ \mbox{and}\ z_k=y_k,\\
  0 & \mbox{for}\ i\geq k, k=1, k=g\ \mbox{or}\ z_k=x_k.
  \end{array}
  \right.
$$

\end{thm}

In order to prove Theorem \ref{main theorem}, we need the following two lemmas. Let $\mathsf{H}_\mathbb{Z}$ be all the real harmonic $1$-forms on $C_0$ with integral periods.

\begin{lem}
\label{eta}
On the curve $C_0$, let $\eta$ be a $1$-form on $C_0$ satisfying the conditions 
$$
\left\{
  \begin{array}{l}
\displaystyle d\eta =\sum_k h_{1, k}\wedge h_{2, k},\\
\displaystyle \int_X \eta\wedge\ast\alpha =0\ \mbox{for}\ \mbox{any}\ \mbox{closed}\ 1\mbox{\mbox{-}form}\ \alpha\ \mbox{on}\ X,\\
\iota^{\ast}\eta=\eta,
  \end{array}
  \right.
$$
where 
$\iota$ is the hyperelliptic involution of $C_0$ and
$h_{1, k}, h_{2, k} \in \mathsf{H}_\Z$ such that 
$\sum_k (h_{1, k}, h_{2, k})
=\sum_k \int_{C_0} h_{1, k}\wedge h_{2, k}=0$.\\
Then for any $j$
$$\int_{e_j}\eta=0.$$
\end{lem}

\begin{proof}
We will have $\eta$ explicitly.
 For any $\displaystyle \sum_k h_{1, k}\wedge h_{2, k}$,
there exist $a_{i, j}^1, a_{i, j}^2\in \C$ such that
$\displaystyle\sum_k h_{1, k}\wedge h_{2, k} =
\sum_{i, j}a_{i, j}^1\omega_i \wedge \overline{\omega_j}+a_{i, j}^2\overline{\omega_i} \wedge \omega_j$,
 where $i,j\in \{1,2, \ldots, g\}$.
The $(1, 1)$-form $\omega_i \wedge \overline{\omega}_j$ is
$\displaystyle{{\lambda^{i-1}\overline{\lambda}^{j-1}}\over{\mu\overline{\mu}}}d\lambda\wedge d\overline{\lambda}$
 in a coordinate $\lambda$ satisfying $\mu^2=\lambda^{2g+2}-1$.
Take a polynomial $f(\lambda, \overline{\lambda})$ of degree at most
 $2g-2$ which belongs to $\C[\lambda, \overline{\lambda}]$
 so that $\displaystyle{{f(\lambda)}\over{\mu\overline{\mu}}}
d\lambda\wedge d\overline{\lambda}=\sum_k h_{1, k}\wedge h_{2, k}$.
 It is clear that $\displaystyle{{f(\lambda)}\over{\mu\overline{\mu}}}d\lambda\wedge d\overline{\lambda}$ is invariant under the action of
 the hyperelliptic involution 
$\iota:(\lambda,\mu) \mapsto (\lambda,-\mu)$,
 since $\mu\overline{\mu}=|\mu^2|=|\lambda^{2g+2}-1|=|(-\mu)^2|$.
So we regard 
$\displaystyle{{f(\lambda)}\over{\mu\overline{\mu}}}
d\lambda\wedge d\overline{\lambda}$ as a $1$-form on $\C P^1$.
On the curve $C_0$, Harris (\cite{H-1} in Section $5, 6$ and \cite{H-2}) gave $\eta$ in the following explicit forms.
\begin{align*}
\eta&={{-1}\over{2\pi}}\int_{\lambda\in\C P^1}\Im{\biggl(}{{dz}\over{z-\lambda}}{\biggr)}\displaystyle{{f(\lambda)}\over{|\lambda^{2g+2}-1|}}d\lambda\wedge d\overline{\lambda}\\
=&{{-1}\over{2\pi}}{{1}\over{2\sqrt{-1}}}{\biggl(}dz
\int{{1}\over{z-\lambda}}\displaystyle{{f(\lambda)}\over{|\lambda^{2g+2}-1|}}d\lambda\wedge d\overline{\lambda}-d\overline{z}\int\overline{{\biggl(}{{1}\over{z-\lambda}}{\biggr)}}\displaystyle{{f(\lambda)}\over{|\lambda^{2g+2}-1|}}d\lambda\wedge d\overline{\lambda}{\biggr)},
\end{align*}
in a coordinate $z$ satisfying $w^2=z^{2g+2}-1$.
It satisfies 
$$
\iota^{\ast}\eta=\eta.
$$
This equation allows us to have
$$\int_{e_k}\eta=\int_{\iota(e_k)^{-1}}\eta
=-\int_{\iota(e_k)}\eta=-\int_{e_k}\iota^{\ast}\eta=-\int_{e_k}\eta.$$
Then we obtain $\int_{e_j}\eta=0$.
\end{proof}

\begin{lem}
\label{volume}
On the curve $C_0$,
$$
\mbox{\phantom{gi}}I(z_i\otimes z_j \otimes z_k) =
\left.
  \begin{array}{cc}
  0 & \mbox{for}\ i\neq j\neq k\neq i,
  \end{array}
\right.
$$
$$
I(x_i\otimes y_i\otimes z_k -x_{k+1}\otimes y_{k+1}\otimes z_k)=
\left\{
  \begin{array}{cc}
  {\displaystyle {1/{2}}} & \mbox{for}\ i<k, k=2, 3, \ldots, g-1\ \mbox{and}\ z_k=y_k,\\
  0 & \mbox{for}\ i\geq k+2, k=1, k=g\ \mbox{or}\ z_k=x_k.
  \end{array}
  \right.
$$
\end{lem}
\begin{proof}
It is enough to consider the iterated integral part of the harmonic volume by Lemma \ref{eta}.

Type $(1)$\\
Lemma \ref{a hyperelliptic iterated integral} gives us $I(z_i\otimes z_j\otimes z_k)\equiv 0$  for $i\neq j\neq k\neq i$.

Type $(2)$\\
We compute $I(x_i\otimes y_i\otimes z_k -x_{k+1}\otimes y_{k+1}\otimes z_k)$ for $i\neq k$ and $i\neq k+1$. When $i<k$, $k=2, 3, \ldots, g-1$ and $z_k=y_k$,
\begin{align*}
&I(x_i\otimes y_i\otimes y_k -x_{k+1}\otimes y_{k+1}\otimes y_k)\\
&=\int_{a_i}\beta_i\beta_k -\int_{a_{k+1}}\beta_{k+1}\beta_{k} \\
&={1\over{2(g+1)^2}}\biggl{\{}-(g+1)\sum_{u=1}^kt_{2i-2u+2}\biggr{\}}-{1\over{2(g+1)^2}}\biggl{\{}-(g+1)\sum_{u=1}^kt_{2(k+1)-2u+2}\biggr{\}}\\
&={1\over{2(g+1)^2}}\{-(g+1)(g-k+1)\}-{1\over{2(g+1)^2}}\{-(g+1)(-k)\}\\
&=-1/{2}\\
&= 1/{2}\ \mathrm{mod}\ \Z.
\end{align*}
It is similarly shown that $I(x_i\otimes y_i\otimes z_k -x_{k+1}\otimes y_{k+1}\otimes z_k)=0$ for $i\geq k+2$, $k=1$, $k=g$ or $z_k=x_k$.
\end{proof}

Before the proof of Theorem \ref{main theorem}, we recall some results about the moduli space of compact Riemann surfaces. Let $\Sigma_g$ be a closed oriented surface of genus $g$. Its mapping class group, denoted here by $\Gamma_g$, is the group of isotopy classes of orientation preserving diffeomorphisms of $\Sigma_g$.
This group acts on the Teichm\"uller space $\mathcal{T}_g$ of $\Sigma_g$ and the quotient space $\mathcal{M}_g$ is the moduli space of Riemann surfaces of genus $g$. The group $\Gamma_g$ acts naturally on the first homology group $H_1(\Sigma_g; \Z)$ of $\Sigma_g$.
 Let $T_g$ be the subgroup of $\Gamma_g$, which acts trivially on $H_1(\Sigma_g; \Z)$ and we call it the Torelli group. Its action on $\mathcal{T}_g$ is free and 
the quotient $\mathcal{I}_g={T_g}\backslash{\mathcal{T}_g}$, 
called the Torelli space, is the moduli space of compact Riemann surfaces 
with a fixed symplectic basis of $H_1(\Sigma_g; \Z)$. There is a natural 
projection $p_T: \mathcal{I}_g\to \mathcal{M}_g.$

 Let $\mathcal{H}_g\subset \mathcal{M}_g$ be the moduli space of hyperelliptic curves of genus $g$.
The hyperelliptic mapping class group $\Delta_g$ is the subgroup of $\Gamma_g$ defined by
$$\{\varphi\in \Gamma_g; \varphi\iota=\iota\varphi\},$$
where $\iota$ is the hyperelliptic involution of $\Sigma_g$.
We choose $\widetilde{\mathcal{H}}_g$ a connected component of $p_T^{-1}(\mathcal{H}_g)$ with the symplectic basis in Figure \ref{figure}.
$\widetilde{\mathcal{H}}_g$ is a complex submanifold of dimension $2g-1$ of $\mathcal{I}_g$.
Let $T_g^H$ denote the group $T_g\cap \Delta_g$.
The moduli space $\mathcal{H}_g$ is known to be connected and has a natural structure of a quasi-projective orbifold.
Hence we have ${\mathcal{H}}_g=p_T(\widetilde{\mathcal{H}}_g)$. 
The group $\Delta_g$ can be considered as its orbifold fundamental group and $T_g^H$ is the fundamental group of $\widetilde{\mathcal{H}}_g$.

\begin{proof}
(Theorem \ref{main theorem})\\
One of the key points of this proof is that the harmonic volume of $C$
 belongs to $\Hom_{\Delta_g}((H^{\otimes 3})^\prime, \Z/{2\Z})=\Hom_{\Z}((H^{\otimes 3})^\prime, \Z/{2\Z})^{\Delta_g}$.
Let $E\to \mathcal{H}_g$ be a flat vector bundle with a fiber $\Hom_{\Z}((H^{\otimes 3})^\prime, \Z/{2\Z})$ and $(p_T|_{\widetilde{\mathcal{H}}_g})^\ast E$ the pullback of the flat vector bundle $E$.
 Harris \cite{H-1} proved that $I$ varies in $\mathcal{I}_g$ continuously.
For any hyperelliptic curves, $I \equiv 0$ or $I \equiv 1/{2}$ modulo $\Z$.
 Hence the flat vector bundle $(p_T|_{\widetilde{\mathcal{H}}_g})^\ast E$
 has a locally constant section $\widetilde{I}$ associated to $I$.
Moreover, $\widetilde{\mathcal{H}}_g$ is arcwise connected and
the monodromy representation
$T_g^H\to \mathrm{Aut}(\Hom_{\Z}((H^{\otimes 3})^\prime, \Z/{2\Z}))$
is trivial.
Therefore $\widetilde{I}$ is constant on $\widetilde{\mathcal{H}}_g$.
Since ${\mathcal{H}}_g=p_T(\widetilde{\mathcal{H}}_g)$, 
the harmonic volumes of hyperelliptic curves can be reduced to
the calculation of $C_0$.
The result follows from Lemma \ref{volume}. 
\end{proof}

\section{The harmonic volumes of hyperelliptic curves \\from a topological viewpoint}
In this section, we study $\Hom_{\Delta_g}((H^{\otimes 3})^{\prime},\Z_2)$ which contains the harmonic volume $I$.
Let $\Z_2$ denote the field $\Z/{2\Z}$.

Birman and Hilden proved the following theorem.
\begin{thm}\label{B-H} $($\cite{B-H}, $\mathrm{Theorem}$ $8)$
The hyperelliptic mapping class group $\Delta_g$ admits the following presentation;
\begin{itemize}
  \item generators: $\sigma_1, \sigma_2, \ldots, \sigma_{2g+1}$
  \item relations:
\begin{description}
    \item[(1)] $\sigma_n\sigma_m=\sigma_m\sigma_n, |n-m|\geq 2,$
    \item[(2)] $\sigma_n\sigma_{n+1}\sigma_n=\sigma_{n+1}\sigma_n\sigma_{n+1}, 1\leq n\leq 2g, $
    \item[(3)] $\theta^{2g+2}=1, $
    \item[(4)] $(\theta\kappa)^2=1, $
    \item[(5)] $\sigma_1(\theta\kappa)=(\theta\kappa)\sigma_1, $
\end{description}
where $\theta=\sigma_1\sigma_2\cdots\sigma_{2g+1}$ and $\kappa=\sigma_{2g+1}\sigma_{2g}\cdots\sigma_1.$
\end{itemize}
\end{thm}
\begin{rem}
The generator $\sigma_i,\ 1\leq i\leq 2g+1,$ is equal to the Dehn twist along
 the simple closed curve $l_i$ in $C$ in Figure \ref{figure4}.
\end{rem}
\begin{center}
\begin{figure}[htbp]
\unitlength 0.1in
\begin{picture}(56.00,16.00)(6.80,-22.00)
%
\special{pn 20}%
\special{pa 1800 600}%
\special{pa 5000 600}%
\special{fp}%
\special{pa 5000 2200}%
\special{pa 1800 2200}%
\special{fp}%
%
\special{pn 20}%
\special{ar 1800 1400 800 800  1.5707963 4.7123890}%
%
\special{pn 20}%
\special{ar 5010 1400 800 800  4.7123890 6.2831853}%
\special{ar 5010 1400 800 800  0.0000000 1.5707963}%
%
\special{pn 20}%
\special{ar 1960 1360 280 160  6.2831853 6.2831853}%
\special{ar 1960 1360 280 160  0.0000000 3.1415927}%
%
\special{pn 20}%
\special{ar 1960 1400 248 104  3.1415927 6.2831853}%
%
\special{pn 20}%
\special{ar 3630 1360 280 160  6.2831853 6.2831853}%
\special{ar 3630 1360 280 160  0.0000000 3.1415927}%
%
\special{pn 20}%
\special{ar 4840 1360 280 160  6.2831853 6.2831853}%
\special{ar 4840 1360 280 160  0.0000000 3.1415927}%
%
\special{pn 20}%
\special{ar 3630 1400 248 104  3.1415927 6.2831853}%
%
\special{pn 20}%
\special{ar 4840 1400 248 104  3.1415927 6.2831853}%
%
\special{pn 8}%
\special{sh 0.600}%
\special{ar 1040 1400 24 24  0.0000000 6.2831853}%
%
\special{pn 8}%
\special{sh 0.600}%
\special{ar 5760 1400 24 24  0.0000000 6.2831853}%
%
\special{pn 8}%
\special{sh 0.600}%
\special{ar 5144 1400 24 24  0.0000000 6.2831853}%
%
\special{pn 8}%
\special{sh 0.600}%
\special{ar 3930 1400 24 24  0.0000000 6.2831853}%
%
\special{pn 8}%
\special{sh 0.600}%
\special{ar 2264 1400 24 24  0.0000000 6.2831853}%
%
\special{pn 8}%
\special{sh 0.600}%
\special{ar 1656 1400 24 24  0.0000000 6.2831853}%
%
\special{pn 8}%
\special{sh 0.600}%
\special{ar 3330 1400 24 24  0.0000000 6.2831853}%
%
\special{pn 8}%
\special{sh 0.600}%
\special{ar 4536 1400 24 24  0.0000000 6.2831853}%
%
\special{pn 20}%
\special{sh 1}%
\special{ar 2820 1400 10 10 0  6.28318530717959E+0000}%
\special{sh 1}%
\special{ar 2740 1400 10 10 0  6.28318530717959E+0000}%
\special{sh 1}%
\special{ar 2900 1400 10 10 0  6.28318530717959E+0000}%
%
\special{pn 8}%
\special{pa 1000 1400}%
\special{pa 680 1400}%
\special{fp}%
\special{pa 5800 1400}%
\special{pa 6280 1400}%
\special{fp}%
%
\special{pn 8}%
\special{ar 6040 1400 160 480  3.7310362 6.2831853}%
\special{ar 6040 1400 160 480  0.0000000 2.3561945}%
%
\special{pn 8}%
\special{pa 5928 1720}%
\special{pa 5912 1680}%
\special{fp}%
\special{sh 1}%
\special{pa 5912 1680}%
\special{pa 5918 1749}%
\special{pa 5932 1730}%
\special{pa 5955 1734}%
\special{pa 5912 1680}%
\special{fp}%
\put(58.4000,-8.8000){\makebox(0,0)[lb]{$180^\circ$}}%
\put(60.4000,-20.4000){\makebox(0,0)[lb]{$\iota$}}%
\put(58.0000,-22.0000){\makebox(0,0)[lb]{$C$}}%
%
\special{pn 13}%
\special{ar 3630 1390 500 370  0.0000000 6.2831853}%
%
\special{pn 13}%
\special{ar 4850 1400 500 370  0.0000000 6.2831853}%
%
\special{pn 13}%
\special{ar 1960 1400 500 370  0.0000000 6.2831853}%
%
\special{pn 13}%
\special{ar 1350 1400 310 240  3.1415927 6.2831853}%
%
\special{pn 13}%
\special{pa 5760 1400}%
\special{pa 5757 1432}%
\special{pa 5749 1463}%
\special{pa 5736 1492}%
\special{pa 5719 1519}%
\special{pa 5699 1544}%
\special{pa 5675 1565}%
\special{pa 5649 1583}%
\special{pa 5621 1600}%
\special{pa 5592 1613}%
\special{pa 5562 1624}%
\special{pa 5531 1631}%
\special{pa 5499 1637}%
\special{pa 5468 1639}%
\special{pa 5436 1640}%
\special{pa 5404 1637}%
\special{pa 5372 1632}%
\special{pa 5341 1625}%
\special{pa 5311 1615}%
\special{pa 5281 1602}%
\special{pa 5254 1585}%
\special{pa 5227 1567}%
\special{pa 5204 1546}%
\special{pa 5183 1521}%
\special{pa 5165 1495}%
\special{pa 5152 1466}%
\special{pa 5143 1435}%
\special{pa 5140 1403}%
\special{pa 5140 1400}%
\special{sp -0.045}%
%
\special{pn 13}%
\special{ar 5450 1400 310 240  3.1415927 6.2831853}%
%
\special{pn 13}%
\special{pa 4550 1400}%
\special{pa 4547 1432}%
\special{pa 4539 1463}%
\special{pa 4526 1492}%
\special{pa 4509 1519}%
\special{pa 4489 1544}%
\special{pa 4465 1565}%
\special{pa 4439 1583}%
\special{pa 4411 1600}%
\special{pa 4382 1613}%
\special{pa 4352 1624}%
\special{pa 4321 1631}%
\special{pa 4289 1637}%
\special{pa 4258 1639}%
\special{pa 4226 1640}%
\special{pa 4194 1637}%
\special{pa 4162 1632}%
\special{pa 4131 1625}%
\special{pa 4101 1615}%
\special{pa 4071 1602}%
\special{pa 4044 1585}%
\special{pa 4017 1567}%
\special{pa 3994 1546}%
\special{pa 3973 1521}%
\special{pa 3955 1495}%
\special{pa 3942 1466}%
\special{pa 3933 1435}%
\special{pa 3930 1403}%
\special{pa 3930 1400}%
\special{sp -0.045}%
%
\special{pn 13}%
\special{ar 4240 1400 310 240  3.1415927 6.2831853}%
%
\special{pn 13}%
\special{pa 1660 1400}%
\special{pa 1657 1432}%
\special{pa 1649 1463}%
\special{pa 1636 1492}%
\special{pa 1619 1519}%
\special{pa 1599 1544}%
\special{pa 1575 1565}%
\special{pa 1549 1583}%
\special{pa 1521 1600}%
\special{pa 1492 1613}%
\special{pa 1462 1624}%
\special{pa 1431 1631}%
\special{pa 1399 1637}%
\special{pa 1368 1639}%
\special{pa 1336 1640}%
\special{pa 1304 1637}%
\special{pa 1272 1632}%
\special{pa 1241 1625}%
\special{pa 1211 1615}%
\special{pa 1181 1602}%
\special{pa 1154 1585}%
\special{pa 1127 1567}%
\special{pa 1104 1546}%
\special{pa 1083 1521}%
\special{pa 1065 1495}%
\special{pa 1052 1466}%
\special{pa 1043 1435}%
\special{pa 1040 1403}%
\special{pa 1040 1400}%
\special{sp -0.045}%
\put(54.0000,-10.0000){\makebox(0,0)[lb]{$l_1$}}%
\put(48.0000,-10.0000){\makebox(0,0)[lb]{$l_2$}}%
\put(42.0000,-10.0000){\makebox(0,0)[lb]{$l_3$}}%
\put(36.0000,-10.0000){\makebox(0,0)[lb]{$l_4$}}%
\put(20.0000,-10.0000){\makebox(0,0)[lb]{$l_{2g}$}}%
\put(14.0000,-10.0000){\makebox(0,0)[lb]{$l_{2g+1}$}}%
\end{picture}%
\caption{}
\label{figure4}
\end{figure}
\end{center}
Let $H_{\Z_2}$ denote $H_1(C; \Z_2)$.
A homomorphism $\rho: \Delta_g\to \mathrm{Sp}(2g; \Z_2)$
 is given by the action on the homology group $H_{\Z_2}$.
 So $H_{\Z_2}$ is a $\Z_2\Delta_g$-module, 
where $\Z_2\Delta_g$ is the group ring of $\Delta_g$.
 We consider $e_i,a_j$ and $b_j$ for $0\leq i\leq 2g+1$
and $1\leq j\leq g$ in Section $3.1$.
The first homology classes of $a_j$ and $b_j$ are
denoted by $x_j$ and $y_j\in H_{\Z_2}$ respectively.
Let $B$ denote the branch locus $\{p_i\}_{i=0,1,\ldots, 2g+1}$.
We deform $e_i$, denoted by $e_i^{\prime}$, to avoid $P_i$
 in a sufficiently small neighborhood of $P_i$
so that $\pi(e_i^\prime)$ surrounds $p_i$ and
$\{\pi(e_i^{\prime})\}_{i=0,1,\ldots, 2g+1}$
 is a generator of $H_1(\C P^1-B;\Z_2)$.
Since the coefficients are in $\Z_2$, the homology class of
$e_i^\prime$ is independent of the choice of $e_i^\prime$.
See Figure \ref{figure5}.

\begin{center}
\begin{figure}[htbp]
\unitlength 0.1in
\begin{picture}( 20.3000,  5.2000)(  7.7000,-12.3000)
%
\special{pn 8}%
\special{sh 0.600}%
\special{ar 1200 1000 46 46  0.0000000 6.2831853}%
%
\special{pn 8}%
\special{sh 0.600}%
\special{ar 2600 1000 46 46  0.0000000 6.2831853}%
\put(7.7000,-12.4000){\makebox(0,0)[lb]{$\pi(Q_0)=\pi(Q_1)$}}%
\put(25.3000,-11.8000){\makebox(0,0)[lb]{$p_i$}}%
\put(16.8000,-8.8000){\makebox(0,0)[lb]{$\pi(e^\prime_i)$}}%
\put(28.0000,-14.0000){\makebox(0,0)[lb]{$\C P^1$}}%
%
\special{pn 13}%
\special{ar 2600 1000 200 200  3.2904826 6.2831853}%
\special{ar 2600 1000 200 200  0.0000000 2.9927027}%
%
\special{pn 13}%
\special{pa 1200 1020}%
\special{pa 2400 1020}%
\special{fp}%
\special{pa 1200 980}%
\special{pa 2400 980}%
\special{fp}%
%
\special{pn 20}%
\special{pa 1800 1020}%
\special{pa 1810 1020}%
\special{fp}%
\special{sh 1}%
\special{pa 1810 1020}%
\special{pa 1744 1000}%
\special{pa 1758 1020}%
\special{pa 1744 1040}%
\special{pa 1810 1020}%
\special{fp}%
\special{pa 1960 980}%
\special{pa 1950 980}%
\special{fp}%
\special{sh 1}%
\special{pa 1950 980}%
\special{pa 2018 1000}%
\special{pa 2004 980}%
\special{pa 2018 960}%
\special{pa 1950 980}%
\special{fp}%
%
\special{pn 20}%
\special{pa 2580 800}%
\special{pa 2570 800}%
\special{fp}%
\special{sh 1}%
\special{pa 2570 800}%
\special{pa 2638 820}%
\special{pa 2624 800}%
\special{pa 2638 780}%
\special{pa 2570 800}%
\special{fp}%
\end{picture}%
\caption{}
\label{figure5}
\end{figure}
\end{center}

 Arnol'd \cite{A} proved the following.
A linear map $\nu:H_{\Z_2}\to H_1(\C P^1-B;\Z_2)$
 defined by $\nu(x_i)=\pi(e_{2i-1}^{\prime})+\pi(e_{2i}^{\prime}),
 \nu(y_i)=\pi(e_0^{\prime})+\pi(e_1^{\prime})+\cdots+\pi(e_{2i-1}^{\prime})$
 is injective.
This map gives the short exact sequence
$$
\xymatrix{
0 \ar[r] & H_{\Z_2} \ar[r]^{\hspace{-35pt}\nu} & H_1(\C P^1-B;\Z_2)
\ar[r] & \Z_2 \ar[r] & 0.
}
$$
Here the map $H_1(\C P^1-B;\Z_2)\to \Z_2$ is the
augmentation map $\pi(e_i^\prime)\mapsto 1$.
Let $f_i$ denote $\pi(e_{0}^{\prime})+\pi(e_{i}^{\prime})$ for 
$i=1,2,\ldots, 2g+1$.
Using $\nu$, we identify $H_{\Z_2}$ with the subgroup of 
$H_1(\C P^1-B;\Z_2)$
generated by $f_1,f_2,\ldots, f_{2g+1}$.
It is clear that $f_1+f_2+\cdots +f_{2g+1}=0$.
A surjective homomorphism $\mu:\Delta_g\to S_{2g+2}$ is defined by 
$\mu(\sigma_j)=(j-1,j)$.
Let $\rho^\prime: S_{2g+2}\to \mathrm{Sp}(2g; \Z_2)$ be the homomorphism
induced by the action on $H_1(\C P^1-B;\Z_2)$ given by the permuting
$\pi(e_{0}^{\prime}),\pi(e_{1}^{\prime}),\ldots,\pi(e_{2g+1}^{\prime})$.
 Arnol'd \cite{A} obtained the commutative diagram
$$
\xymatrix{
H_{\Z_2} \ar[rr]^{\rho(\sigma_j)} \ar[d]^{\nu} & 
                                        & H_{\Z_2} \ar[d]^{\nu} \\
H_1(\C P^1-B;\Z_2)\hspace{10pt}  \ar[rr]^{\rho^\prime(j-1,j)} &
                  &\hspace{10pt}  H_1(\C P^1-B;\Z_2).
}
$$
We identify the actions of $\sigma_1,\sigma_2
\ldots, \sigma_{2g}$ and $\sigma_{2g+1}$
 on $H_{\Z_2}$ with those of the transpositions 
$(0,1),(1,2),\ldots,(2g-1,2g)$ and $(2g,2g+1)$ on 
$H_1(\C P^1-B; \Z_2)$ respectively.

We denote by $\Delta_g^\prime=\{\sigma\in \Delta_g;\ \sigma(P_0)=P_0\}$ and
$\Delta_g^{\prime\prime}=\{\sigma\in \Delta_g;\ \sigma(P_0)=P_0\ \mathrm{and}\ \sigma(P_1)=P_1\}$.
We have $\mu(\Delta_g^{\prime})=S_{2g+1}$ and $\mu(\Delta_g^{\prime\prime})
=S_{2g}$, 
where $S_{2g+1}=\{\sigma\in S_{2g+2};\sigma(\pi(e_0^\prime))=\pi(e_0^\prime)\}$
and $S_{2g}=\{\sigma\in S_{2g+2};\sigma(\pi(e_0^\prime))=\pi(e_0^\prime)
\ \mathrm{and}\ \sigma(\pi(e_1^\prime))=\pi(e_1^\prime)\}$.
As in the proof of Theorem \ref{main theorem},
 the pointed harmonic volume $I_{P_0}$ is an element of $\Hom_{\Delta_g^\prime}(K\otimes H, \Z_2)$.
For a $\Z\Delta_g^{\prime}$-module $M$, we denote
$M^\ast=\Hom_{\Z}(M, \Z_2)$, which is naturally regarded as a
$\Z_2\Delta_g^{\prime}$-module.
Clearly we have $H^\ast=H_{\Z_2}^\ast$.

The homomorphism of short exact sequences
$$
\xymatrix{
0 \ar[r]& K\otimes H \ar[r]& H^{\otimes 3} \ar[r]^{\ (\ ,\ )\otimes \id} & H \ar[r] & 0\\
0 \ar[r]& (H^{\otimes 3})^\prime \ar[r] \ar[u] & H^{\otimes 3} \ar[r]^p \ar[u] & H^{\oplus 3} \ar[r] \ar[u] & 0,
}$$
induces the homomorphism of long exact sequences,
\begin{multline}\label{commutative diagram}
{\scriptsize
\hspace{-10pt}\xymatrix{
H^0(S_{2g+1}; H^\ast) \ar[r]& H^0(S_{2g+1}; (H^{\otimes 3})^\ast) \ar[r] & H^0(S_{2g+1}; (K\otimes H)^\ast) \ar[r] & H^1(S_{2g+1}; H^\ast) \\
H^0(S_{2g+1}; (H^{\oplus 3})^\ast) \ar[r] \ar[u] & H^0(S_{2g+1}; (H^{\otimes 3})^\ast) \ar[r] \ar[u] & H^0(S_{2g+1}; ((H^{\otimes 3})^\prime)^\ast)  \ar[u] \ar[r] & H^1(S_{2g+1}; (H^{\oplus 3})^\ast). \ar[u]
}}
\end{multline}

\begin{lem}
\label{basis representaion}
We have
$$H^0(S_{2g+1}; H^\ast)=0.$$
\end{lem}
\begin{proof}
We take $\varphi\in H^0(S_{2g+1}; H^\ast)$. 
Since $\varphi$ is $S_{2g+1}$-equivariant, $\varphi(f_1)=\varphi(f_2)=
\cdots =\varphi(f_{2g+1})$.
Using $f_1+f_2+\cdots +f_{2g+1}=0$, we have 
 $0=\varphi(f_1+f_2+\cdots +f_{2g+1})=(2g+1)\varphi(f_1)=\varphi(f_1)$.
From $\varphi(f_i)=0,\ 1\leq i\leq 2g+1$, $H^0(S_{2g+1}; H^\ast)=0$ follows.
\end{proof}

We recall the notion of induced and co-induced modules.
Let $\Ind_{S_{2g}}^{S_{2g+1}}\Z_2$
 denote the induced module $\Z_{2}S_{2g+1}\otimes_{\Z_2S_{2g}}\Z_2$
 and $\Coind_{S_{2g}}^{S_{2g+1}}\Z_2$
 the co-induced module $\Hom_{S_{2g}}(\Z_2S_{2g+1}, \Z_2)$.
They are $(2g+1)$-dimensional vector spaces over $\Z_2$.
We denote by $r_i=(i, 1)\otimes 1\in\Ind_{S_{2g}}^{S_{2g+1}}\Z_2$
for $i=1,2,\ldots, 2g+1$.
Then $\{r_i\}_{i=1,2,\ldots, 2g+1}$
 is a basis of $\Ind_{S_{2g}}^{S_{2g+1}}\Z_2$.
 Since $[S_{2g+1}: S_{2g}]<\infty$, we have a natural isomorphism
$\lambda:\Coind_{S_{2g}}^{S_{2g+1}}\Z_2\to\Ind_{S_{2g}}^{S_{2g+1}}\Z_2$
given by $\lambda(s)=\sum_{i=1}^{i=2g+1}(i,1)\otimes s((i,1))$
for $s\in\Coind_{S_{2g}}^{S_{2g+1}}\Z_2$.
Let $s_i$ be the element of $\Coind_{S_{2g}}^{S_{2g+1}}\Z_2$
such that $\lambda(s_i)=r_i$.
We have a natural exact sequence
\begin{equation}\label{exact sequence}
\xymatrix{0 \ar[r] & H_{\Z_2} \ar[r]^{\hspace{-25pt}\phi} & \Coind_{S_{2g}}^{S_{2g+1}}\Z_2 \ar[r]^{\hspace{25pt}\chi} & \Z_2 \ar[r] & 0},
\end{equation}
where $\phi(n_2f_2+n_3f_3+\cdots +n_{2g+1}f_{2g+1})
=n_2s_1+n_3s_2+\cdots+n_{2g+1}s_{2g}+(n_2+n_3+\cdots+n_{2g+1})s_{2g+1}$
and $\chi$ is the augmentation map.

A transfer map is defined as follows.
The canonical surjection $\tau$ of $S_{2g+1}$-modules
$\Ind_{S_{2g}}^{S_{2g+1}}\Z_2\to \Z_2$ is defined by
$\tau(\sigma\otimes a)=\sigma a=a$.
By Shapiro's lemma, we obtain
$H^i(S_{2g+1}; \Coind_{S_{2g}}^{S_{2g+1}}\Z_2)= H^i(S_{2g}; \Z_2)\ $
 for any $i$.
A transfer map
cor$_{S_{2g}}^{S_{2g+1}}: H^i(S_{2g}; \Z_2)\to H^i(S_{2g+1}; \Z_2)$
 is induced by Shapiro's lemma and the following composite mapping
$$\xymatrix{
 \Coind_{S_{2g}}^{S_{2g+1}}\Z_2 \ar[r]^{\hspace{5pt}\lambda} &
 \Ind_{S_{2g}}^{S_{2g+1}}\Z_2 \ar[r]^{\hspace{15pt}\tau} & \Z_2.
}$$
It immediately follows that
$\chi$ is equal to $\tau\circ\lambda$.

\begin{lem}\label{extension}
We have
$$H^1(S_{2g+1}; H^\ast)=0.$$
\end{lem}
\begin{proof}
The exact sequence (\ref{exact sequence}) induces the exact sequence
$$\xymatrix{0 \ar[r] & H^0(S_{2g+1}; H_{\Z_2}) \ar[r]^{\hspace{-25pt}{\phi^\ast}} & H^0(S_{2g+1}; \Coind_{S_{2g}}^{S_{2g+1}}\Z_2) \ar[r]^{\hspace{25pt}{\chi^\ast}} & H^0(S_{2g+1}; \Z_2)}$$
$$\xymatrix{\mbox{\phantom{g}}\ar[r] & H^1(S_{2g+1}; H_{\Z_2}) \ar[r]^{\hspace{-25pt}{\phi^\ast}} & H^1(S_{2g+1}; \Coind_{S_{2g}}^{S_{2g+1}}\Z_2) \ar[r]^{\hspace{25pt}{\chi^\ast}} &  H^1(S_{2g+1}; \Z_2).}$$
 By Shapiro's lemma, we obtain
$H^i(S_{2g+1}; \Coind_{S_{2g}}^{S_{2g+1}}\Z_2)= H^i(S_{2g}; \Z_2)\ $
 for $i=0, 1$.
We have $H^0(S_{2g+1}; \Z_2)=\Z_2$ and $H^0(S_{2g}; \Z_2)=\Z_2$,
since the actions of $S_{2g+1}$ and $S_{2g}$ on $\Z_2$ are trivial.
Let $\mathrm{sign}_i$ be the signature map $S_{i}\to \Z_2$ for $i=2g,2g+1$.
Since $2g, 2g+1\geq 6> 5$,
we obtain $H^1(S_{i}; \Z_2)=\Z_2$ and $\mathrm{sign}_i$
generates $H^1(S_{i}; \Z_2)$ for $i=2g,2g+1$.
In order to prove $H^1(S_{2g+1}; H^\ast)=H^1(S_{2g+1}; H_{\Z_2})=0$,
it is enough to prove that $\chi^\ast$ is an isomorphism.
Let $1_{i}$ denote the nontrivial element of $H^0(S_{i}; \Z_2)$
for $i=2g,2g+1$.
Since $\chi=\tau\circ\lambda$, we have
 $\chi^\ast=\mathrm{cor}_{S_{2g}}^{S_{2g+1}}: H^0(S_{2g}; \Z_2)\to
 H^0(S_{2g+1}; \Z_2)$.
Lemma \ref{basis representaion} gives 
$H^0(S_{2g+1};H_{\Z_2})=H^0(S_{2g+1};H^{\ast})=0$.
Then we obtain $\mathrm{cor}_{S_{2g}}^{S_{2g+1}}(1_{2g})=1_{2g+1}$
and the isomorphism $\chi^\ast: H^0(S_{2g}; \Z_2)\to H^0(S_{2g+1}; \Z_2)$.
We apply the transfer formula
\begin{center}
cor$_{S_{2g}}^{S_{2g+1}}(\mathrm{res}_{S_{2g}}^{S_{2g+1}}(\mathrm{sign}_{2g})
\cup 1_{2g})=\mathrm{sign}_{2g}\cup\mathrm{cor}_{S_{2g}}^{S_{2g+1}}(1_{2g})$
\end{center}
to $\mathrm{sign}_{2g}$ and $1_{2g}$.
So cor$_{S_{2g}}^{S_{2g+1}}\mathrm{res}_{S_{2g}}^{S_{2g+1}}(\mathrm{sign}_{2g})
=\mathrm{sign}_{2g}$.
Since $\chi^\ast$ is surjective, we have the isomorphism
 $\chi^\ast=\mathrm{cor}_{S_{2g}}^{S_{2g+1}}: H^1(S_{2g}; \Z_2)\to
 H^1(S_{2g+1}; \Z_2)$.
 Then $H^1(S_{2g+1}; H^\ast)=H^1(S_{2g+1}; H_{\Z_2})=0$.
\end{proof}

Using the diagram (\ref{commutative diagram}), Lemma \ref{basis representaion} and Lemma \ref{extension}, we get the homomorphism of the commutative diagram
$$\xymatrix{
 H^0(S_{2g+1}; (H^{\otimes 3})^\ast)
                \ar[r] & H^0(S_{2g+1}; (K\otimes H)^\ast)  \\
 H^0(S_{2g+1}; (H^{\otimes 3})^\ast) \ar[r] \ar[u]
        & H^0(S_{2g+1}; ((H^{\otimes 3})^\prime)^\ast)). \ar[u] 
}$$
The two horizontal and one left-hand vertical homomorphisms are isomorphisms.
Then the other right-hand vertical homomorphism is an isomorphism.
We have $H^0(S_{2g+1}; (H^{\otimes 3})^\ast)
=H^0(S_{2g+1}; ((H^{\otimes 3})^\prime)^\ast))
=H^0(S_{2g+1}; (K\otimes H)^\ast)$.

\begin{lem}\label{hom}
$$H^0(S_{2g+1}; (H^{\otimes 3})^\ast)=\Z_2.$$
Moreover, the unique nontrivial element $\psi\in H^0(S_{2g+1}; (H^{\otimes 3})^\ast)$ is an $S_{2g+1}$-homomorphism $H^{\otimes 3}\to \Z_2$ defined by
$$\psi(f_i\otimes f_j\otimes f_k)=
\left\{
\begin{array}{ll}
0 &for\ i\neq j\neq k\neq i,\\
0 &for\ i=j=k,\\
1 &otherwise.
\end{array}
\right.$$
\end{lem}

\begin{proof}
Let $\psi$ be an element of $H^0(S_{2g+1}; (H^{\otimes 3})^\ast)$.
Since $\psi$ is $S_{2g+1}$-equivariant,
there exist $a, b_1, b_2, b_3$ and $c\in \Z_2$ such that
$\psi(f_i\otimes f_i\otimes f_i)=a$,
$\psi(f_j\otimes f_i\otimes f_i)=b_1,
\psi(f_i\otimes f_j\otimes f_i)=b_2, 
\psi(f_i\otimes f_i\otimes f_j)=b_3$ for $i\neq j$ and 
$\psi(f_i\otimes f_j\otimes f_k)=c$ for $i\neq j\neq k\neq i$.
The dimension of $\Z_2$-vector space $H^0(S_{2g+1}; (H^{\otimes 3})^\ast)$
is not greater than $1$.
Since $I(x_i\otimes x_i\otimes y_i -x_i\otimes
x_{i+1}\otimes y_{i+1} -x_{i+1}\otimes x_i\otimes y_{i+1})\equiv 1/{2}$,
 $I_{P_0}\in H^0(\Delta_g^\prime; ((H^{\otimes 3})^\prime)^\ast))
=H^0(S_{2g+1}; ((H^{\otimes 3})^\ast)$ is not $0$.
We obtain $H^0(S_{2g+1}; ((H^{\otimes 3})^\prime)^\ast))\neq 0$.
Hence we have $H^0(S_{2g+1}; (H^{\otimes 3})^\ast)=\Z_2$.
It is clear that the generator of $H^0(S_{2g+1}; (H^{\otimes 3})^\ast)$
 is $\psi$ as above.
\end{proof}

\begin{cor} \label{0-hom}
$$H^0(S_{2g+2}; (H^{\otimes 3})^\ast)=0.$$
\end{cor}
\begin{proof}
Take $\psi\in H^0(S_{2g+1}; (H^{\otimes 3})^\ast)$ in the proof of Lemma \ref{hom}.
Let $b$  denote $b_1=b_2=b_3$.
Using $\rho(\sigma_1)(f_i)=f_1+f_i$ for $i=2,3,\ldots, 2g+1$, we have
 $0=\psi(f_2\otimes f_3\otimes f_4)=\psi(\rho(\sigma_1)
(f_2\otimes f_3\otimes f_4))=3b=b$.
The equation $a=b_1=b_2=b_3=c=0$ gives $H^0(S_{2g+2}; (H^{\otimes 3})^\ast)=0$.
\end{proof}

Using the diagram (\ref{commutative diagram}), Lemma \ref{basis representaion},
 Lemma \ref{extension} and Lemma \ref{hom}, we have
\begin{prop}\label{p.h.v.}
$$H^0(\Delta_g^\prime; (K\otimes H)^\ast)=H^0(\Delta_g^\prime; ((H^{\otimes 3})^\prime)^\ast)=\Z_2.$$
\end{prop}

This gives us the following theorem.
\begin{thm}\label{homological main theorem}\ 
$$H^0(\Delta_g; ((H^{\otimes 3})^\prime)^\ast)= \Z_2.$$
\end{thm}
\begin{proof}
We have a natural injection 
$H^0(\Delta_g; ((H^{\otimes 3})^\prime)^\ast)\hookrightarrow
 H^0(\Delta_g^\prime; ((H^{\otimes 3})^\prime)^\ast)$.
Using Proposition \ref{p.h.v.}, the dimension of $\Z_2$-vector space 
$H^0(\Delta_g; ((H^{\otimes 3})^\prime)^\ast)$ is not greater than $1$.
As in the proof of Lemma \ref{hom},
the harmonic volume $I\in H^0(\Delta_g; ((H^{\otimes 3})^\prime)^\ast)$
is not $0$.
Hence $H^0(\Delta_g; ((H^{\otimes 3})^\prime)^\ast)= \Z_2$.
\end{proof}

\begin{proof}
(The second proof of Theorem \ref{main theorem})\\
Using Theorem \ref{homological main theorem}, Proposition \ref{p.h.v.}
and Lemma \ref{hom},
we identify $H^0(S_{2g+1}; (H^{\otimes 3})^\ast)$ with 
$H^0(\Delta_g; ((H^{\otimes 3})^\prime)^\ast)$, whose generator is regarded
as $\psi$ in Lemma \ref{hom}.
We substitute 
$$
\left\{\begin{array}{l}
\displaystyle x_i=f_{2i-1}+ f_{2i} ,\\
\displaystyle y_i=f_1+f_2+\cdots +f_{2i-1},
\end{array}\right.
$$
for elements of the type $(1)$ and $(2)$ in Section $2$.
Then the direct computation of $\psi$ gives us Theorem \ref{main theorem}.
\end{proof}

 The harmonic volume $I$ gives a geometric interpretation of a theorem established by Tanaka.

\begin{thm}
\label{Tanaka}
$(\mathrm{Tanaka}$\cite{T}, $\mathrm{Theorem}$ $1.1)$\\
If $g\geq 2$, then
$$H_1(\Delta_g; H)=\Z_2.$$
\end{thm}

Tanaka obtained the generator of $H_1(\Delta_g; H)$, using the relations of $\Delta_g$ in Theorem \ref{B-H}.
 Since $\Delta_g$ acts transitively on $H$, $H_0(\Delta_g, H)=0$.
By the universal coefficient theorem,
$$H^1(\Delta_g; H^\ast)=\Hom_{\Z}(H_1(\Delta_g; H); \Z_2)=\Z_2.$$
We have $H^1(\Delta_g; H^\ast)=\Z_2$.

 By Corollary \ref{0-hom}, it is clear that $H^0(\Delta_g; (H^{\otimes 3})^\ast)=0$. The short exact sequence
 $$\xymatrix{0 \ar[r] & (H^{\otimes 3})^\prime \ar[r] & H^{\otimes 3} \ar[r]^p & H^{\oplus 3}\ar[r] & 0}$$
 gives us a connected homomorphism $\delta: H^0(\Delta_g; ((H^{\otimes 3})^\prime)^\ast)\to H^1(\Delta_g; (H^{\oplus 3})^\ast)$ and it is injective.
Since $I$ is $S_3$-invariant,
we may consider $\delta I=(\delta I|_H,\delta I|_H,\delta I|_H)
\in H^1(\Delta_g; H^\ast)^{\oplus 3}$.
Here $\delta I|_H$ is the restriction
$H^0(\Delta_g; ((H^{\otimes 3})^\prime)^\ast)\to
H^1(\Delta_g; H^\ast)$.

\begin{prop}
The generator of $H^1(\Delta_g; H^\ast)$ is $\delta I|_H$.
\end{prop}
\begin{proof}
If $\delta I|_H$ is {\it not} the generator of $H^1(\Delta_g; H^\ast)$,
we have $\delta I=0\in H^1(\Delta_g; (H^{\oplus 3})^\ast)$.
This contradicts that $\delta$ is injective.
\end{proof}

\end{document}